\documentclass[12pt]{article}

\usepackage{amsmath,amsfonts,amssymb}

\newtheorem{theorem}{Theorem}[section]
\newtheorem{lemma}{Lemma}[section]
\newtheorem{remark}{Remark}[section]
\newtheorem{corollary}{Corollary}[section]

\def\no{\noindent}
\def\pa{\partial}

\def\Box{\diamond}

\def\ds{\displaystyle}

\def\Om{\Omega}  
\def\om{\omega}
\def\ga{\gamma}
\def\eps{\varepsilon}

\newcommand{\R}{\mathbb R}

\begin{document}

\centerline{\Large\bf The Fourier Singular Complement Method}

\medskip
\centerline{\Large\bf for the Poisson problem. Part I: prismatic domains}

\vskip1truecm

\centerline{P. Ciarlet, Jr,
\footnote{
ENSTA \& CNRS UMR 2706, 32, boulevard Victor, 75739 Paris Cedex 15, France. 
{\em This author was supported in part by France/Hong Kong Joint Research Scheme}.
}
\quad 
B. Jung,
\footnote{Department of Mathematics, Chemnitz University of Technology, 
D-09107 Chemnitz, Germany.
{\em This author was supported by DGA/DSP-ENSTA 00.60.075.00.470.75.01 Research Programme}.
}
\quad 
S. Kaddouri,
\footnote{
ENSTA \& CNRS UMR 2706, 32, boulevard Victor, 75739 Paris Cedex 15, France. 
}
\quad 
S. Labrunie,
\footnote{
IECN, Universit\'e Henri Poincar\'e Nancy~I \& INRIA (Projet CALVI), 54506 
Vand\oe uvre-l\`es-Nancy cedex, France. 
}
\quad 
J. \ Zou
\footnote{Department of Mathematics, The Chinese University of Hong Kong, 
Shatin, N.T., Hong Kong. 
{\em This author was fully supported by Hong Kong RGC grants (Project CUHK4048/02P and
project 403403)}.
}
}
   
\abstract{\footnotesize This is the first part of a threefold article, aimed at solving 
numerically the Poisson problem in three-dimensional prismatic or axisymmetric domains.
In this first part, the Fourier Singular Complement Method is introduced and analysed, in 
prismatic domains. In the second part, the FSCM is studied in axisymmetric domains with conical 
vertices, whereas, in the third part, implementation issues, numerical tests and comparisons with 
other methods are carried out. The method is based on a Fourier expansion in the direction 
parallel to the reentrant edges of the domain, and on an improved variant of the Singular 
Complement Method in the 2D section perpendicular to those edges. 
Neither refinements near the reentrant edges of the domain nor cut-off functions are required in 
the computations to achieve an optimal convergence order in terms of the mesh size and the number 
of Fourier modes used.}

\bigskip
{\it Date of this version} :  April 18, 2005 \\

\section{Introduction}

 The {\em Singular Complement Method} (SCM) was originally introduced by Assous {\em et al} 
\cite{AsCS98,AsCS00}, for the 2D static or instationary Maxwell equations without charges. 
The cases with charges have been recently solved by Garcia {\em et al} \cite{AsCG00,Garc02}, 
including the numerical solution to the 2D Vlasov-Maxwell system of equations.
The SCM has been extended in \cite{CiHe03} to the 2D Poisson problem. Further extensions to 
the 2D heat or wave equations, or to similar problems with piecewise constant coefficients, 
can be obtained easily. As a matter of fact, this stems from the analysis which is performed 
hereafter (see Remark~\ref{rmk-wave}).
 The primary basis of the SCM is the decomposition of the solution into regular and singular 
parts. Methodologically speaking, the SCM consists in adding some singular test functions to
the usual $P_1$ Lagrange FEM so that one recovers the optimal $H^1$-convergence rate, even 
in non-convex domains. In 2D, one may simply add one singular test function per reentrant 
corner. \\

 There exist a couple of numerical methods in the literature for accurately solving 2D Poisson 
problems in non-convex domains. It was shown in \cite{CiHe03} that the SCM can be reformulated 
so that it coincides with the approach of Moussaoui \cite{Mous84,AmMo89} when L-shaped domains are
considered. The SCM differs from the {\em Dual Singular Function Method} (DSFM) of Blum and 
Dobrowolski \cite{BlDo82} in that it requires no cut-off functions. Actually, when the numerical 
implementation of the SCM is carried out, the cut-off function is traded for a non-homogeneous 
boundary condition. Note that Cai and Kim \cite{CaKi01} recently proposed a new SFM which 
involves the evaluations of singular and cut-off functions and the solution of a nonsymmetric
elliptic problem. The SCM is clearly different from (anisotropic) mesh refinement techniques
\cite{Raug78,ApHe94,Hein96,ApNi98,Apel99}, and can be applied efficiently to instationary 
problems (see Remark~\ref{rmk-wave}), since it does not need mesh refinement and thus larger 
timesteps may be allowed. However the anisotropic mesh refinement methods have one advantage: 
they require only a partial knowledge of the most singular part of the solution. \\

 The numerical solution of 3D singular Poisson problems is quite different from the 2D case, 
and much more difficult. This is a relatively new field of research: most approaches rely on 
anisotropic mesh refinement, see for instance \cite{ApHe94,ApSW96,Hein96,ApNi98,HeNW00}, and
\cite{Apel99} and Refs. therein. 
To our knowledge, this series of papers is the first attempt to generalize the SCM for 
three-dimensional singular Poisson problems. Specifically, we shall consider the numerical 
solution of the Poisson problem: \\
{\em Find $u\in H^1_0(\Om)$ such that}
 \begin{equation}\label{Poisson}
-\Delta u =f \quad \mbox{in} \quad \Omega, 
 \end{equation}
where $f\in L^2(\Omega)$, and 
$\Om$ is a (right) {\em prismatic domain} described by 
\begin{equation}
 \Om =\om\times Z\,, \label{prismatic}
\end{equation} 
and $\om$ is a two-dimensional general polygonal domain, $Z$ is an interval varying from $0$ to 
a positive constant $L$ on the $x_3$-axis. The {\em bases} of the domain are the subsets of the 
boundary $\pa\Om$, which are included in the planes $\{x_3=0\}$ and $\{x_3=L\}$. \\

The case of an {\em axisymmetric domain} is considered in the companion paper \cite{CJKL+04b}. 
When the Poisson problem (\ref{Poisson}) is solved in this class of domains, two difficulties 
arise. The first difficulty is that one has to deal with weighted Sobolev spaces, the weights 
being functions of the distance to the axis. The second one is that there exist two kinds of 
{\em geometrical singularities}: reentrant edges like in the prismatic case, and, in addition, 
{\em sharp conical vertices}. As for {\em implementation issues and comparisons} with other 
methods (such as variants of our method, the FSCM, or mesh refinement techniques \cite{HeNW00}), 
we refer the reader to~\cite{CJKL+04c}. \\

 The rest of the paper is organized as follows. In the next Section, some theoretical results 
concerning the regularity of the solution to the Poisson problem in prismatic domains are recalled. 
{\em A priori} regularity results of the solution $u$ to (\ref{Poisson}), and a first splitting of 
the solution into regular and singular parts, are emphasized. In Section 3, some results 
about the Fourier expansion along $x_3$ are recalled and/or proven. This suggests a framework for 
building the {\em Fourier Singular Complement Method} (FSCM) for accurately solving the problem 
(\ref{Poisson}), using a Fourier expansion in $x_3$, and an improved variant of the Singular
Complement Method \cite{CiHe03} in the 2D section $\omega$. In Section~\ref{error-estimate}, we study 
the variant of the SCM, based on a theoretical splitting of the solution $u_\mu$ to 2D problems of the
form $-\Delta u_\mu + \mu\,u_\mu = f_\mu$ in $\om$ (with a parameter $\mu\ge0$ related to the Fourier 
modes). The main feature of the regular-singular splitting is that it is chosen {\em independently of 
$\mu$}; this independence is important, and very helpful, from the computational point of view.
Estimates on Sobolev norms of $u_\mu$ and its splitting are established. To end this Section, 
the SCM is considered from a numerical point of view, to approximate $u_\mu$ accurately, {\em via} 
the discretization of the splitting: the optimal $H^1$-norm convergence of the order $O(h)$ is 
recovered. In the last Section, we first prove a refined splitting of the solution $u$ to the 
3D Poisson problem under suitable assumptions on the right-hand side $f$, using the Fourier 
expansion along $x_3$. Then, we build the numerical algorithms which define the FSCM, and we show that 
the FSCM has the optimal convergence of order $O(h+N^{-1})$, where $h$ is the 2D mesh size and $N$ is 
the number of Fourier modes used. \\

Throughout this paper, when two quantities $a$ and $b$ are such that $a\le C\, b$, with a constant $C>0$
which depends only on the geometry of the domain, we shall use the notation $a\lesssim b$.

\section{Poisson problem in prismatic domains}
 Let us recall that a (right, open) cylinder of $\R^3$, with axis parallel to $x_3$, is equal to $D\times 
I$, where $D$ is any connected (open) subset of $\R^2$, and $I$ is any (open) interval of $\R$.
 Let us proceed then with some remarks on the class of domains $\Om$, i.e., the prismatic domains. 
{\em A priori}, such domains could be considered:
 \begin{itemize}
     \item either as truncated infinite cylinders;
     \item or as polyhedra.
  \end{itemize}
 \indent As it happens, considering $\Om$ as a {\em polyhedron} is helpful, in a simple manner. 
Indeed, from \cite{ApNi98,CoDa00}, we know that, in any polyhedra, the solution $u$ to (\ref{Poisson})
can be split~as
 \begin{eqnarray}
     &&u=u_r+u_e+u_v, \mbox{ with }u_r\in H^2(\Om),\label{general-splitting}\\
     &&u_e=\sum_e\mu_e(\rho_e,z_e)\sin(\alpha_e\phi_e),\mbox{ and } 
       u_v=\sum_v\sum_{-1/2<\lambda_v<1/2}\mu_{v,\lambda_v}\rho_v^{\lambda_v}\Phi_v(\theta_v,\phi_v)\nonumber.
 \end{eqnarray}
 Above, $u_r$ is called the {\em regular} part, $u_e$ the {\em edge singularity} part, and $u_v$ 
the {\em vertex singularity} part. Note that when $u_e\neq0$ or $u_v\neq0$, they do not belong to 
$H^2(\Om)$. The summation in $u_e$ is taken over all {\em reentrant} edges $e$, $(\rho_e,\phi_e,z_e)$ denote 
the local cylindrical coordinates, and $\pi/\alpha_e$ the dihedral angle (so that 
$\alpha_e\in]1/2,1[$). 
Last, the summation in $u_v$ is taken over all {\em non-convex} vertices $v$ and over all eigenvalues 
$\lambda_v$ of the Laplace-Beltrami operator, which belong to the interval $]-1/2,1/2[$, and 
$(\rho_v,\theta_v,\phi_v)$ denote the local spherical coordinates. \\
 In our case, i.e., when $\Om$ is a {\em prismatic} domain with polygonal bases, 
it has been shown \cite{StWh88,Apel99} that the vertex singularity part $u_v$ always vanishes, so 
(\ref{general-splitting}) reduces to
 \begin{equation}\label{prismatic-splitting}
     u=u_r+u_s, \mbox{ with }u_r\in H^2(\Om)\mbox{ and }u_s=\sum_e\mu_e(\rho_e,z_e)\sin(\alpha_e\phi_e).
 \end{equation}
 \indent Let us describe how one can fall into the other class, that of the {\em infinite cylinders}. \\
 The first step is to introduce a suitable continuation $\tilde u$ of the solution $u$ (odd reflection 
at the bases) along the $x_3$ direction from $Z$ to $\R$: one builds a problem to be solved in the 
infinite cylinder $C_\infty=\om\times\R$. Unfortunately, with this continuation technique, one gets 
a solution (and data) which is not in $L^2(C_\infty)$.
 Thus, one introduces in a second step a smooth truncation function $\eta$, such that $\eta(x_3)$ is 
equal to one for $x_3\in]-L/2,3L/2[$, and to zero for $|x_3|>2L$. Then, one multiplies $\tilde u$ 
by $\eta$, to obtain a Poisson problem in $C_\infty$ with solution $u^\eta=u\,\eta$. This time, one has 
$u^\eta\in H^1(C_\infty)$ (and $f^\eta=-\Delta u^\eta\in L^2(C_\infty)$).
 By construction, the restriction of $u^\eta$ on $\Om$ coincides with $u$.

 Interestingly, it has been proven in \cite{Gris87,Lenc93}, that a splitting similar to 
(\ref{prismatic-splitting}) holds for $u^\eta$. Furthermore, $u^\eta_s$ can be expressed as 
 \begin{equation} \label{stress}
 u^\eta_s = \ga^\eta_e(\rho_e, x_3)\rho_e^{\alpha_e}\sin(\alpha_e\phi_e).
\end{equation} 
The function $\ga^\eta_e$ in (\ref{stress}) is often called in mechanics the {\em stress intensity 
distribution}. 
On the one hand, in the original paper \cite{Gris87}, $\ga^\eta_e$ is expressed as a convolution product. 
On the other hand, in \cite{Lenc93}, it is characterized as the solution to a second order PDE. 
Finally, the regularity of the singular part $u^\eta_s$, can be expressed accurately as follows
\cite{Lenc93}. Let $\delta$ denote the minimal distance between two reentrant edges, and for each 
reentrant edge $e$, let $\Om_e=\{\vec x\in\Om\ :\ d(\vec x,e)<\delta/2\}$. Then
 \begin{equation}\label{reg-u-2}
 \left\{
 \begin{array}{l}
 u\in H^{1+\alpha-\eps}(\Om),\ \forall \eps>0,\ \alpha=\min_e\alpha_e,\\
 u\in H^2(\Om\setminus\cup_e\bar\Om_e),\\
 \rho_e^{\beta_e}\pa_i u\in H^1(\Om_e),\ \forall e,\ \forall\beta_e>1-\alpha_e,\ i=1,2,\\
 \pa_3 u\in H^1(\Om).
 \end{array}
 \right.
 \end{equation}
 \indent In Section \ref{FSCM}, using a Fourier expansion along the $x_3$ axis, we recover some properties
which are very similar to (\ref{prismatic-splitting}-\ref{reg-u-2}). \\

 We end this Section with remarks on other possible boundary conditions. \\
 
 If the boundary condition for $u$ on the bases of the physical domain $\Om$ are the non-homogeneous 
Dirichlet boundary condition:
\[
  u=g \quad \mbox{at} \quad x_3=0 \quad \mbox{and} \quad x_3=L, 
\]
one can set $w=u-\tilde g$ with $\tilde g$ being a continuation of $g$ into $\Om$. Then the problem 
reduces to the case with the solution $w$ satisfying the homogeneous Dirichlet boundary condition, 
assuming that $\tilde g\in H^2(\Om)$. \\

If the boundary condition for $u$ is the homogeneous Neumann boundary condition $\pa_n u=0$ on 
$\pa\Om$, then one can replace the $\sin(\alpha_e\phi_e)$ factor in (\ref{prismatic-splitting}) by the 
expected $\cos(\alpha_e\phi_e)$. Moreover, to obtain an expression like (\ref{stress}), one uses an 
even reflection of $u$ at the bases of the domain. 
If we have the non-homogeneous Neumann boundary condition: 
\[
  \pa_n u = g \quad \mbox{at} \quad x_3=0 \quad \mbox{and} \quad x_3=L, 
\]
one may then study the solution
\[
  w(\vec x) = u(\vec x)-\int_0^{x_3} \tilde g(x_1, x_2,z)\,dz 
\]
first, which satisfies the homogeneous Neumann boundary conditions at the bases of the domain $\Om$. 
Here $\tilde g$ is a continuation of $g$ into $\Om$ which is again assumed to belong to $H^2(\Om)$.  \\

 From now on, we assume, for ease of exposition, that the polygon $\om$ has only one reentrant corner 
$C$, i.e., with an interior angle larger than $\pi$, denoted as $\pi/\alpha$, with $1/2<\alpha<1$. In 
particular, the summation which defines the singular part $u_s$ in (\ref{prismatic-splitting}) 
reduces to exactly one term.

\section{Fourier expansion}
We devote this Section to some justifications about the Fourier series expansion of the 
Poisson solution to (\ref{Poisson}). First, one can show, following for instance Heinrich's 
proof of Lemma~3.2 in~\cite{Hein93}, the well-known result
\begin{lemma}\label{lem:f1}
For any $f\in L^2(\Omega)$, there exist Fourier coefficients defined by 
\begin{equation} \label{f1}
 f_k(x_1, x_2) = 
  \frac{2}{L} \int_0^{L} f(x_1,x_2,x_3) \sin \frac{k\pi}{L} x_3\,dx_3,\quad k=1,2,3,\cdots, 
\end{equation}
such that $f_k\in L^2(\om)$ and 
\begin{equation} \label{f2}
 f(x_1, x_2, x_3) =\sum_{k=1}^\infty f_k(x_1, x_2) \sin \frac{k\pi}{L}x_3
 \quad\mbox{ a.e. in }\Omega \,,
\end{equation}
and 
\begin{equation} \label{f3}
 \|f\|^2_{L^2(\Omega)} = \frac{L}2 \sum_{k=1}^\infty \|f_k\|^2_{L^2(\om)} < \infty. 
\end{equation}
If $f\in H^1_0(\Omega)$, then $f_k\in H^1_0(\omega)$ for all $k$ and 
\begin{equation}\label{f-H1_norm}
 \|\nabla f\|^2_{L^2(\Omega)} =\frac{L}2 \sum_{k=1}^\infty 
   \Big\{ \|\nabla f_k\|^2_{L^2(\om)} + \Big( \frac{k\pi}{L}\Big)^2 \|f_k\|^2_{L^2(\om)}
   \Big\} <\infty.
\end{equation}
\end{lemma}
 For $f$ in $L^2(\Om)$, let us introduce the sequence of partial sums $(F_K)_K$ of the Fourier 
decomposition of $f$, which converges to $f$ in $L^2(\Om)$, cf. (\ref{f2}): \\ 
\begin{equation}\label{partial-sum}
    F_K=\ds\sum_{k=1}^K f_k\sin \frac{k\pi}{L}x_3,\mbox{ for }K>0.
\end{equation}
 We note that when $f$ is in $H^1_0(\Om)$, $(F_K)_K$  converges to $f$ in $H^1_0(\Om)$, according 
to (\ref{f-H1_norm}). Also, the sine functions can be replaced by cosine functions with the same 
argument $\ds\frac{k\pi}{L}x_3$, and (\ref{f1}-\ref{f-H1_norm}) still holds (for (\ref{f-H1_norm}), 
with any $f$ in $H^1(\Om)$.) \\

 In our subsequent analysis, summations like 
 \begin{equation}\label{sum-k4}
     \sum_{k=1}^\infty k^4 \|f_k\|^2_{L^2(\om)}
 \end{equation}
 will appear. The result below provides a characterization of elements $f$ of $L^2(\Om)$, which 
are such that (\ref{sum-k4}) is bounded. Let the following Sobolev spaces be introduced:
 \begin{eqnarray*}
  h^1(\Om)&:=&H^1(]0,L[,L^2(\om))
            = \{f\in L^2(\Om)\ :\ \pa_3 f\in L^2(\Om)\}~;\\
  h_{\diamond}^1(\Om)
          &:=&H^1_0(]0,L[,L^2(\om))
            = \{f\in h^1(\Om)\ :\ f_{|\{x_3=0\}}=f_{|\{x_3=L\}}=0\}~;\\
  h^2(\Om)&:=&H^2(]0,L[,L^2(\om))
            = \{f\in h^1(\Om)\ :\ \pa_{33}f\in L^2(\Om)\}.
     \end{eqnarray*}
\begin{lemma}\label{lem:f2} Given $f\in L^2(\Om)$, one has the following equivalences
 \begin{eqnarray}
     &&f\in h_{\diamond}^1(\Om) \Longleftrightarrow 
       \sum_{k=1}^\infty k^2 \|f_k\|^2_{L^2(\om)}<\infty~;\label{h10}\\
     &&f\in h_{\diamond}^1(\Om)\cap h^2(\Om)\Longleftrightarrow 
       \sum_{k=1}^\infty k^4 \|f_k\|^2_{L^2(\om)}<\infty.\label{h10-cap-h2}
 \end{eqnarray}
\end{lemma}
  
 \noindent {\it Proof}. Let $f$ be in $L^2(\Om)$. \\
 Assume in addition that $f\in h_{\diamond}^1(\Om)$. We note that, by the definition of the 
Fourier mode $f_k$ and integration by parts ($f$ vanishes at the bases), one has
 \[ (k\pi)f_k = -2 \int_0^{L} f\left(\cos \frac{k\pi}{L} x_3\right)'\,dx_3
              =  2 \int_0^{L} \pa_3 f\cos \frac{k\pi}{L} x_3\,dx_3.\] 
 Since by assumption, $\pa_3 f$ is in $L^2(\Om)$, one gets the expected
 $\ds\sum_{k=1}^\infty k^2 \|f_k\|^2_{L^2(\om)}<\infty$. \\
 Let us prove the reciprocal assertion. For $f$ in $L^2(\Om)$, as the sequence $(F_K)_K$ (see 
(\ref{partial-sum})) converges to $f$ in $L^2(\Om)$, one infers that $(\pa_3 F_K)_K$ converges 
to $\pa_3 f$ in $H^{-1}(\Om)$. Now, if the sum is bounded, $\ds\left(\pa_3 F_K\right)_K$ 
is a Cauchy sequence in $L^2(\Om)$, so it converges in this space, and its limit $\pa_3 f$ 
is in $L^2(\Om)$. Since both $(F_K)_K$ and $(\pa_3 F_K)_K$ converge in $L^2(\Om)$, $(F_K)_K$ 
converges to $f$ in $h^1(\Om)$, and, as $F_K$ belongs to $h_{\diamond}^1(\Om)$ for all $K$, 
$f$ is also in $h_{\diamond}^1(\Om)$, which proves (\ref{h10}). \\

 In order to establish (\ref{h10-cap-h2}), one proceeds similarly, by performing a second 
integration by parts. Note that for this additional integration by parts, no assumption is required 
on the trace of $f$ at the bases, since the $\ds(\sin\frac{k\pi}{L}x_3)_k$ vanish there.
  \[ \frac{k^2\pi^2}{L}f_k 
 = 2 \frac{k\pi}{L}\int_0^{L} \pa_3 f\cos \frac{k\pi}{L} x_3\,dx_3
 = -2 \int_0^{L} \pa_{33} f\sin \frac{k\pi}{L} x_3\,dx_3 .\] 
 With this identity, one concludes the proof easily. \hfill $\Box$ \\
 
Now that the general results have been obtained, we focus on the Poisson 
problem~(\ref{Poisson}). Consider the weak form of the Poisson problem:
\begin{equation} \label{Poissonw}
a(u, v)=f(v) \quad \forall\,v\in H_0^1(\Om)
\end{equation} 
where $a(\cdot, \cdot)$ and $f(\cdot)$ are given by 
\[
 a(u, v) =\int_\Om \nabla u \cdot \nabla v\,dx, \quad f(v)=\int_\Om f\,v\,dx. 
\]
We expand the solution $u$ in (\ref{Poisson}) in the Fourier sine series:
\begin{equation}\label{fourier0}
 u(x_1, x_2, x_3)=\sum_{k=1}^\infty u_k(x_1, x_2) \sin \frac{k\pi}{L}x_3.
\end{equation} 
Following again Heinrich's proof of Lemma~3.2 in~\cite{Hein93}, the next two Lemmas hold.
\begin{lemma}\label{lem:f3}
For any $u, v\in H_0^1(\Om)$, we have 
\[
 a(u, v) =\frac{L}2 \sum_{k=1}^\infty a_k(u_k, v_k), \quad 
 f(v) =\frac{L}2 \sum_{k=1}^\infty f_k(v_k), 
\] 
where $a_k$ and $f_k$ are given by 
\[ 
 a_k(u_k, v_k) =\int_\om\Big\{ \nabla u_k \cdot \nabla v_k + 
  \Big( \frac{k\pi}{L}\Big)^2 u_kv_k\Big\} dx_1dx_2, 
 \quad f_k(v_k) =\int_\om f_k\,v_k\,dx_1dx_2, 
\] 
and $u_k$, $v_k$ and $f_k$ are Fourier coefficients of $u$, $v\in H_0^1(\Om)$ 
and $f\in L^2(\Om)$ respectively. 
\end{lemma}

\begin{lemma}\label{lem:f4}
For any $f\in L^2(\Om)$, let $u\in H_0^1(\Om)$ be the unique weak solution 
of (\ref{Poissonw}) and $u_k$ and $f_k$ be the Fourier coefficients of 
$u$ and $f$. Then $u_k\in H_0^1(\om)$ is the unique solution 
of the following 2D weak problem: \\
{\em Find $u_k\in H^1_0(\om)$ such that}
\begin{equation} \label{Poissonk}
 a_k(u_k, v) =f_k(v) \quad \forall\, v\in H_0^1(\om). 
\end{equation}
Moreover, $u_k$ satisfies the following {\em a priori} estimates: 
\begin{eqnarray*} 
 \int_\om\Big\{ |\nabla u_k |^2 + \Big( \frac{k\pi}{L}\Big)^2 u^2_k\Big\} dx_1dx_2 
  &\le& \Big( \frac{L}{k\pi}\Big)^2 \|f_k\|^2_{L^2(\om)}, \quad k=1, 2, \cdots,\\
  \sum_{k=1}^\infty k^2 
  \Big\{ \|\nabla u_k \|^2_{L^2(\om)} + \Big( \frac{k\pi}{L}\Big)^2 
  \|u_k\|^2_{L^2(\om)} \Big\} 
   &\le& \frac{2 L}{\pi^2} \|f\|^2_{L^2(\Om)}. 
\end{eqnarray*} 
\end{lemma}
 This means that the $k$-th Fourier mode of $u$ is characterized as the unique solution to the 2D 
problem\\
{\em Find $u_k\in H^1_0(\om)$ such that}
\begin{equation} \label{2dprob}
 -\Delta u_k +\left( \frac {k\pi}{L} \right)^2 u_k =f_k \quad 
 \mbox{in} \quad \om; \quad u_k=0 \quad \mbox{on} \quad {\pa \om}. 
\end{equation} 

 As Corollaries, one gets a convergence result of the sequence of partial sums $(U_K)_K$ of 
the Fourier decomposition of $u$, and also the last result of (\ref{reg-u-2}). 
\begin{corollary}\label{cor:f1}
Let $f\in L^2(\Om)$, and $u$ be the solution to (\ref{Poisson}). Then $(U_K)_K$ converges to $u$ in 
$H^1(\Om)$, and $(\Delta U_K)_K$ converges to $-f$ in $L^2(\Om)$.
\end{corollary}
 \noindent {\it Proof}. The fact that $(U_K)_K$ converges to $u$ in $H^1(\Om)$ is a consequence
of Lemma \ref{lem:f1}. Then, one notes that
 \[-\Delta U_K=\sum_{k=1}^K 
 (-\Delta u_k+\left( \frac {k\pi}{L} \right)^2 u_k)\sin\frac{k\pi}{L}x_3
 \stackrel{(\ref{2dprob})}{=}\sum_{k=1}^Kf_k\sin\frac{k\pi}{L}x_3=F_K,\]
 which yields the result on the convergence of $(\Delta U_K)_K$. \hfill $\Box$ \\
\begin{corollary}\label{cor:f2}
Let $f\in L^2(\Om)$, and $u$ be the solution to (\ref{Poisson}). Then $\pa_3u\in H^1(\Om)$.
\end{corollary}
 \noindent {\it Proof}. We prove that, for $i=1,2,3$, $\pa_{i3}u$ belongs to $L^2(\Om)$. \\
 For $i=3$, thanks to the last bound of the Lemma~\ref{lem:f4}, there holds
$\ds\sum_{k=1}^\infty k^4 \|u_k\|^2_{L^2(\om)}<\infty$. Result (\ref{h10-cap-h2}) yields 
$u\in h_{\diamond}^1(\Om)\cap h^2(\Om)$, so that $\pa_{33}u$ is in $L^2(\Om)$. \\
 For $i=1,2$, we note that
 \[\pa_{i3}U_K=-\frac{\pi}{L}\sum_{k=1}^K k\pa_iu_k\cos\frac{k\pi}{L}x_3.\]
 According again to the last estimate in Lemma~\ref{lem:f4},
$\ds\sum_{k=1}^\infty k^2 \|\pa_iu_k\|^2_{L^2(\om)}<\infty$, so $\pa_{i3}u$ is in $L^2(\Om)$. 
\hfill $\Box$ \\

 To conclude this Section, we note that the Fourier expansion (\ref{fourier0}) of $u$ together with 
the series of 2D problems (\ref{2dprob}) suggest the numerical approximation scheme below, i.e.,
define the Fourier SCM (FSCM) approximation of the solution $u$ to (\ref{Poissonw}) as follows: 
 \begin{equation} \label{FSCM-approx}
 U_N^h (x_1, x_2, x_3) =\sum_{k=1}^N u_{k}^h(x_1, x_2) 
 \sin \frac{k\pi}{L}x_3\,
 \end{equation}
 where $N$ is the total number of Fourier modes used in the approximation, and $u_{k}^h$ is a 
suitable approximation of $u_k$, to be studied in the next two Sections.

 \section{Regular-singular decomposition in the 2D domain $\om$: theoretical study}
\label{error-estimate}
The main interest of this paper is to propose some efficient
numerical method for solving the three-dimensional singular Poisson
problem (\ref{Poisson}) in a prismatic domain. Basically, the method reduces the
3D problem into a series of 2D Poisson-like problems, see (\ref{2dprob}), 
by the Fourier expansion of the 3D solution along the $x_3$-direction.
This Section will thus focus on the 2D singular Poisson problem: \\
{\em Find $u_\mu\in H^1_0(\om)$ such that}
\begin{equation}\label{laplace1}
-\Delta u_\mu\ + \mu\,u_\mu\ = f\ \quad \mbox{in} ~~\om.
\end{equation}
In the case of the Fourier expansion, one considers $\mu=k^2\pi^2/L^2$ and $f=f_k$ 
in (\ref{laplace1}). Due to the presence of the Fourier mode index $k$, the coefficient 
$\mu$ varies in a large range, from $\pi^2/L^2$ to $N^2\pi^2/L^2$, where $N$ is 
the number of Fourier modes required subsequently in the numerical approximation 
(cf. Section~\ref{FSCM}). This brings in one of the main difficulties in the subsequent
error estimates, which should hold for all $\mu$'s in a large range. \\ 

 As a preliminary remark, we note that, according to \cite{Gris92}, 
the most singular part of the solution to (\ref{laplace1}) is of the form 
$\rho^\alpha\sin(\alpha\theta)$, compared to (\ref{prismatic-splitting}) in 3D. \\

 Let $\ga_1$, $\ga_2$, $\cdots$, $\ga_K$  be the line segments of 
$\pa\om$, where $\ga_1$ and $\ga_2$ are two line segments 
which form the single re-entrant corner of $\om$.
Our numerical method is based on the following 
important decomposition of the space $L^2(\om)$ \cite{Gris92}:
\begin{equation}\label{decomp10}
 L^2(\om)=\Delta [H^2(\om)\cap H_0^1(\om)] \stackrel\perp\oplus N\,,
\end{equation}
where $N$ is a space of singular harmonic functions defined by 
\[
N=\Big\{p\in L^2(\om)\ :\ \Delta p=0,\ 
p|_{\ga_k}=0 \mbox{ in } (H_{00}^{1/2}(\ga_k))',\ 1\le k\le K\Big\}.
\]
 Above, the space $H^{1/2}_{00}(\ga_k)$ is made up of elements of $H^{1/2}(\ga_k)$, such that 
their continuation to $\pa\om$ by zero belongs to $H^{1/2}(\pa\om)$. Its dual space is denoted 
by $(H_{00}^{1/2}(\ga_k))'$. To understand that the boundary condition on $p$ holds in this dual 
space, let us mention that one can prove that, given any $\tilde\phi$ in $H^2(\om)\cap 
H_0^1(\om)$, $\pa_n\tilde\phi|_{\ga_k}$ belongs to $H^{1/2}_{00}(\ga_k)$. Then, the fact that 
$p|_{\ga_k}=0$ simply reflects a surjectivity property, which states that the mapping 
$\tilde\phi\mapsto\pa_n\tilde\phi|_{\ga_k}$ is onto, from $H^2(\om)\cap H_0^1(\om)$ to 
$H^{1/2}_{00}(\ga_k)$. \\

As the domain $\om$ has only one re-entrant corner, we know dim$(N)=1$, and $N$=span$\{p_s\}$ for 
some $p_s\in N\setminus\{0\}$, see Grisvard \cite{Gris92}.

Let $\phi_s$ be an element in $H^1_0(\om)$, which solves 
the Poisson problem 
\begin{equation} \label{2ps}
  - \Delta\phi_s =p_s \quad \mbox{in} \quad \om\,.
\end{equation}
Then by the decomposition (\ref{decomp10}), we 
can split  the solution $u_\mu$ to equation (\ref{laplace1}) as
\begin{equation}\label{observ1}
u_\mu=\tilde u_\mu + c_\mu\phi_s,
\end{equation}
where 
$\tilde u_\mu\in H^2(\om)\cap H_0^1(\om)$, and is called 
the regular part of $u_\mu$. 

We will devote the rest of this Section 
to the derivation of some {\em a priori} estimates for the 
solution $u_\mu$, its regular part $\tilde u_\mu$ and the singularity 
coefficient $c_\mu$, as well as the solvability of 
$\tilde u_\mu$ and $c_\mu$. Let us first introduce some notation. 

Throughout the rest of the paper, 
$\alpha_0$ will be a frequently used 
fixed positive constant lying in the interval $]\frac 12, \alpha[$, where 
$\alpha\in ]\frac 12, 1[$ is the singularity exponent. 
$|\cdot|_s$ is used to denote 
the semi-norm of the Sobolev space $H^s(\om)$ for any $s>0$, 
$(\cdot, \cdot)$
and $\|\cdot\|_0$ are used to denote the inner product and the norm 
in the space $L^2(\om)$. Also, $(\cdot, \cdot)$ will be used for 
the dual pairing between the space $H_0^1(\om)$ and $H^{-1}(\om)$ when necessary.

The following lemma summarizes some a priori estimates on $u_\mu$ and $c_\mu$.
\begin{lemma}\label{lem:uk} Let $u_\mu$ be the solution $u_\mu$ to the Poisson problem
(\ref{laplace1}), then we have the following {\em a priori} estimates:
\begin{eqnarray}
\mu\,\|u_\mu\|_0\le \|f\|_0\,, \quad 
\sqrt{\mu}\,|u_\mu|_1 &\le& \frac 1{\sqrt{2}} \|f\|_0\,, \quad 
\|\Delta u_\mu\|_0 \le 2\,\|f\|_0\,, \label{eq:uk1}\\
|c_\mu|&\lesssim&\mu^{-\frac{1-\alpha}{2}} \,\|f\|_0\, \label{asympc}\\
|u_\mu|_{1+\alpha_0} &\lesssim&\mu^{-\frac{1-\alpha_0}2} \,\|f\|_0\,.
\label{asympu}
\end{eqnarray}
\end{lemma}
\noindent {\it Proof}. Multiplying equation (\ref{laplace1}) by $u_\mu$ and integrating over $\om$ 
yield 
\[
  |u_\mu|_1^2 +\mu\,\|u_\mu\|_0^2 \le \|f\|_0\,\|u_\mu\|_0\,, 
\] 
which proves the first estimate in (\ref{eq:uk1}). 
Then applying the Cauchy-Schwarz inequality, we further obtain 
\[
  |u_\mu|_1^2 + \mu\,\|u_\mu\|_0^2 \le \frac 12 \mu\,\|u_\mu\|_0^2 + \frac 1{2\mu} 
  \|f\|_0^2, 
\] 
which leads to the $H^1$ semi-norm estimate in (\ref{eq:uk1}). 

The last estimate in (\ref{eq:uk1}) follows 
immediately from $\Delta u_\mu=\mu u_\mu-f$ and the first inequality 
in (\ref{eq:uk1}). 

\medskip
As far as (\ref{asympc}) is concerned, it is a simple matter to check that the 
singularity coefficient $c_\mu$, multiplied by some constant $\beta^\star$, equals the singularity
coefficient $c(\mu)$ of \cite[pp. 62-69]{Gris92}. 
Indeed, in Grisvard's papers, $u_\mu$ is decomposed into:
 \begin{equation}\label{Grisvard-decomp}
     u_\mu=u_\mu^G+c(\mu){\rm e}^{-\sqrt{\mu}\rho}\xi(\rho)\rho^\alpha\sin(\alpha\phi),\quad u_\mu^G\in 
     H^2(\om)\cap H^1_0(\om)
 \end{equation}
where $\xi$ is a smooth cut-off function, 
equal to one in a neighborhood of 0. 
 
 On the other hand one can decompose the singular part 
in (\ref{observ1}) as (cf. \cite{CiHe03} or 
 (\ref{gamma1}) below)
 \[c_\mu\phi_s=c_\mu\left(\tilde\phi+\beta^\star \rho^\alpha\sin(\alpha\phi)\right),
 \quad \tilde\phi\in H^2(\om), \quad \beta^\star =\frac1\pi \|p_s\|_0^2.\]
Using this, (\ref{observ1}) and (\ref{Grisvard-decomp}), we can write 
 \begin{eqnarray}
&&(c_\mu\beta^\star-c(\mu)\xi(\rho))\rho^\alpha\sin(\alpha\phi)\nonumber\\
&=& u_\mu -(\tilde u_\mu +c_\mu\tilde\phi) - 
c(\mu)\xi(\rho)\rho^\alpha\sin(\alpha\phi)\nonumber \\
 &=&u_\mu^G+c(\mu)\left({\rm e}^{-\sqrt{\mu}\rho}-1\right)
\xi(\rho)\rho^\alpha\sin(\alpha\phi) -(\tilde u_\mu+c_\mu\tilde\phi).
\label{compare}
\end{eqnarray}
Noting that each term on the right-hand side of (\ref{compare}) belongs to $H^2(\om)$, we must have 
 $c_\mu=c(\mu)/\beta^\star$. 
But it is shown in \cite[ineq. (2.5.5)]{Gris92} that \begin{equation}\label{Grisvard-estimate1}
|c(\mu)|\lesssim\mu^{-\frac{1-\alpha}2} \,\|f\|_0,
\end{equation}
which implies (\ref{asympc}). \\
 
 In order to derive the estimate (\ref{asympu}), we shall use
(\ref{Grisvard-decomp}-\ref{Grisvard-estimate2}), with the 
additional norm estimate 
\cite[ineq. (2.5.4)]{Gris92} on the regular part $u_\mu^G$, namely
 \begin{equation}\label{Grisvard-estimate2}
|u_\mu^G|_2+\sqrt\mu|u_\mu^G|_1+\mu\|u_\mu^G\|_0\lesssim\|f\|_0.
\end{equation}
 Indeed, from the estimates 
 \[|u_\mu^G|_1\lesssim\mu^{-1/2}\|f\|_0,\qquad
   |u_\mu^G|_2\lesssim\|f\|_0\,, 
\]
we have then by standard interpolation theory that 
\[
   |u_\mu^G|_{1+\alpha_0}\lesssim\mu^{-\frac{1-\alpha_0}2}\|f\|_0.
\]
Next, we use 
(\ref{Grisvard-estimate1}) and a direct estimate of
the $H^{1+\alpha_0}$ semi-norm to bound the singular part in 
(\ref{Grisvard-decomp}). 
Actually, there holds
 \[ |v|_{1+\alpha_0}^2=\int_{\vec x\in\om}\int_{\vec x'\in\om}
 \frac{|\nabla v(\vec x)-\nabla v(\vec x')|^2}{|\vec x-\vec x'|^{2+2\alpha_0}}
 d\om(\vec x)\,d\om(\vec x'),
 \quad \forall\,v\in H^{1+\alpha_0}(\om).\]
 Due to the uniform smoothness (in $\mu$) of ${\rm e}^{-\sqrt{\mu}\rho}\xi(\rho)\rho^\alpha\sin(\alpha\phi)$
for $\rho\ge \rho_0>0$, it is possible to evaluate the integrals only on
$\om_\infty=\{(\rho,\phi)\in]0,\rho_0[\times]0,\pi/\alpha[\}$.
Then, one performs the changes of variables $s=\sqrt\mu \rho$, $s'=\sqrt\mu \rho'$, to find
 \[|{\rm e}^{-\sqrt{\mu}\rho}\xi(\rho)\rho^\alpha\sin(\alpha\phi)|_{H^{1+\alpha_0}(\om_\infty)}
     \le C(\alpha_0)\mu^{-\frac{\alpha-\alpha_0}{2}}.\]
  This with (\ref{asympc}) leads to 
(\ref{asympu}).  \hfill $\Box$

Now, let us study the solvability of $\tilde u_\mu$ and $c_\mu$ in 
decomposition (\ref{observ1}). 
For convenience, we introduce the notation $a_\mu(\cdot, \cdot)$ and 
the norm $\|\cdot\|_a$:
\[
 a_\mu(w, v) =(\nabla w, \nabla v)+ \mu\, (w, v)\,, \quad
 \|v\|^2_a = a_\mu(v, v)\,,
\]
and the linear mapping $A_\mu$ from $H_0^1(\om)$ to $H^{-1}(\om)$, defined 
by $A_\mu u=-\Delta u +\mu u$, or equivalently by
\[
 _{H^{-1}(\om)}<A_\mu w, v>_{H^1_0(\om)} = a_\mu(w, v) \quad \forall\,w, v\in H_0^1(\om).
\]
It is not difficult to verify that $A_\mu$ is a one-to-one and 
onto mapping, so it is invertible.

So, we claim that $\tilde u_\mu$ and $c_\mu$ solve the following coupled system: 
\begin{eqnarray}
&&a_\mu(\tilde u_\mu, v) +c_\mu\,a_\mu(\phi_s, v)
=(f, v) \quad \forall v\in H^1_0(\om)\,, \label{cmu3}\\
&&\big( \Vert p_s\Vert_0^2 + \mu\vert\phi_s\vert_1^2 \big)\,c_\mu+
\mu\,(\tilde u_\mu, p_s)=( f, p_s)\,. \label{cmu1}
\end{eqnarray}
In fact, by multiplying the equation (\ref{laplace1}) 
by $p_s$ and integrating over $\om$ we obtain  
\[
-(\Delta u_\mu, p_s)+\mu\,(u_\mu, p_s)=( f, p_s)\,,
\]
then (\ref{cmu1}) follows readily from the decomposition (\ref{observ1}), 
the orthogonality between $p_s$ and $\Delta \tilde u_\mu$, 
along with the relation (\ref{2ps}) and its following direct consequence 
\begin{equation}\label{phips}
|\phi_s|_1^2 =(\phi_s, p_s).  
\end{equation}
On the other hand, the solution $u_\mu$ of (\ref{laplace1}) also satisfies the weak form: 
\[ 
(\nabla u_\mu, \nabla v)+\mu\,(u_\mu, v)=(f, v)\,
\quad \forall v\in H^1_0(\om). 
\] 
This and the decomposition (\ref{observ1}) lead to the equation (\ref{cmu3}). 

\medskip
Below, we show the well-posedness of the system (\ref{cmu3})-(\ref{cmu1}). 
\begin{lemma}\label{lem:bound1} 
There exists a unique solution $(\tilde u_\mu, c_\mu)$ to the 
coupled system (\ref{cmu3})-(\ref{cmu1}) and the following stability estimates 
hold: 
\begin{eqnarray*}
 \|\tilde u_\mu\|_a
 &\le& \sqrt2 \,\Big(2\sqrt{\mu} C_P^2+\frac 1{\sqrt{\mu}}\Big) \|f\|_0\,, \\ 
|c_\mu|&\le& 2\,\frac{\|f\|_0}{\|p_s\|_0}\,, \quad \quad
|\tilde u_\mu|_2\le 4\,\|f\|_0\,, 
\end{eqnarray*}
where $C_P$ is the constant in the Poincar\'e inequality.
\end{lemma}

\no {\it Proof}. To see the unique existence, 
we rewrite (\ref{cmu3}) as the following operator form:
\begin{equation}
  A_\mu\tilde u_\mu + c_\mu\,A_\mu\phi_s =f \quad \mbox{in} \quad H^{-1}(\om). \label{oper1}
\end{equation}
As the inverse of $A_\mu$ exists, we know from (\ref{oper1}) that   
$\tilde u_\mu$ can be determined if $c_\mu$ is available: 
\begin{equation}\label{oper4}
   \tilde u_\mu= A_\mu^{-1} f - c_\mu\,\phi_s\,.
\end{equation}
This is exactly our original decomposition (\ref{observ1}). 
Substituting this into (\ref{cmu1}), 
\[
\Big ( \Vert p_s\Vert_0^2 + \mu\vert\phi_s\vert_1^2 \Big )c_\mu+
\mu\,(A_\mu^{-1} f - c_\mu\,\phi_s, p_s)=(f, p_s)\,.
\]
With (\ref{phips}), we obtain that
\begin{equation}
c_\mu =\frac{(f- \mu\,A_\mu^{-1} f, \,p_s)}{\Vert p_s\Vert_0^2}\,. 
\label{oper3}
\end{equation}
With $c_\mu$ uniquely determined, $\tilde u_\mu$ is clearly 
uniquely determined by (\ref{cmu3}) or (\ref{oper4}).

Next, we derive the stability estimates in Lemma~\ref{lem:bound1}. 
We show that these estimates are the consequences of (\ref{oper4}-\ref{oper3}) 
and of the following inequality 
\begin{equation}
\|A_\mu^{-1}\,g\|_0 \le \frac 1\mu\,\|g\|_0 \quad \forall\,g\in 
 L^2(\om).  \label{oper6}
\end{equation} 
In fact, if (\ref{oper6}) is true, then the desired estimate 
on $c_\mu$ follows from (\ref{oper3}): 
\[
 |c_\mu| \le \frac{\|f\|+ \mu\,\|A_\mu^{-1} f\|_0}{\Vert p_s\Vert_0}\le 
 2 \frac{\|f\|_0}{\|p_s\|_0}\,.
\]
On the other hand, we have from (\ref{phips}) and the 
Poincar\'e inequality that 
\[ 
 \|\phi_s\|_0 \le C_P\,\|\nabla \phi_s\|_0\le C_P^2 \|p_s\|_0\,. 
\] 
Using this and the bound of $c_\mu$, we derive 
from (\ref{cmu3}) by taking $v=\tilde u_\mu$ that
\begin{eqnarray*}
 \|\nabla \tilde u_\mu\|_0^2 +\mu \|\tilde u_\mu\|_0^2 
 &\le& \|f\|_0 \|\tilde u_\mu\|_0 + |c_\mu|\, 
 (\|\nabla \phi_s\|_0\|\nabla \tilde u_\mu\|_0+\mu\,\|\phi_s\|_0 
  \|\tilde u_\mu\|_0)\\
 &\le& \|f\|_0 \|\tilde u_\mu\|_0+ 2C_P \|f\|_0\|\nabla \tilde u_\mu\|_0
   + 2\,\mu\,C_P^2 \|f\|_0\|\tilde u_\mu\|_0\,.
\end{eqnarray*}
Then, an application of the Young inequality yields 
\begin{eqnarray*}
\|\nabla \tilde u_\mu\|_0^2 +\mu \|\tilde u_\mu\|_0^2 
 &\le&  \frac 12 \mu\,\|\tilde u_\mu\|_0^2 + \frac 1\mu \|f\|_0^2 + 
    \frac 12 \|\nabla \tilde u_\mu\|_0^2 + 2C_P^2 \|f\|_0^2 
  + 4\mu\,C_P^4 \|f\|_0^2\,. 
\end{eqnarray*}
This implies 
\begin{eqnarray*}
    \frac12\|\tilde u_\mu\|_a^2
    &\le&\left(\frac1\mu+2C_P^2+4\mu C_P^4\right)\|f\|_0^2
    \le\left(\frac1{\sqrt\mu}+2\sqrt\mu C_P^2\right)^2\|f\|_0^2\,,
\end{eqnarray*}
so the desired estimate on $\|\tilde u_\mu\|_a$ follows.

We now show the $H^2$-norm estimate. By the decomposition 
(\ref{oper4}), we have 
$u_\mu=A_\mu^{-1}f=\tilde u_\mu + c_\mu \phi_s$, and 
\[
 -\Delta \tilde u_\mu =-\Delta u_\mu +c_\mu \Delta \phi_s = 
  f-\mu\,u_\mu -c_\mu p_s, 
\] 
which gives 
\[
 \|\Delta \tilde u_\mu\|_0\le \|f\|_0+\mu\|u_\mu\|_0+|c_\mu|\,\|p_s\|_0. 
\]
But we know from Lemma~\ref{lem:uk} that $\mu\|u_\mu\|_0\le \|f\|_0$. This, along with the previous bound
for $c_\mu$, leads to 
\[\|\Delta \tilde u_\mu\|_0\le 4\|f\|_0.\]
Now, for any $\vec v\in H^1(\om)^2$ such that $\vec v\cdot\vec\tau=0$ on $\pa\om$,
with $\vec\tau$ the vector tangential to $\pa\om$, it is well-known (cf. \cite{Cost91}) 
that (since $\om$ is a polygon)
\[\sum_{1\le k,l\le 2}\|\pa_k v_l\|_0^2=\|{\rm curl}\vec v\|_0^2+\|{\rm div}\vec v\|_0^2.\]
So, by taking $\vec v=\nabla\tilde u_\mu$, one actually finds 
\[ 
|\tilde u_\mu|_2=\|\Delta \tilde u_\mu\|_0\le 4\|f\|_0.
\] 

Finally,  it remains to prove (\ref{oper6}). By the definition of $a_\mu(\cdot, \cdot)$, 
we easily see the following lower bound:
\begin{equation}\label{oper7}
 _{H^{-1}(\om)}<A_\mu v, v>_{H^1_0(\om)} =
 a_\mu(v, v) \ge \mu\,\|v\|_0^2 \quad \forall\,v\in H_0^1(\om)\,.
\end{equation}
Then for any $g\in L^2(\om)\subset H^{-1}(\om)$, let $v=A_\mu^{-1}\,g\in H_0^1(\om)$. One has 
$A_\mu\,v=g$ in $L^2(\om)$ and it follows from (\ref{oper7}) that 
\[
 \|A_\mu^{-1}\,g\|_0^2 =\|v\|_0^2 \le \frac 1\mu\,(A_\mu\,v, v)
                       =\frac 1\mu\,(g, A_\mu^{-1}\,g)=\frac 1\mu\|g\|_0\,\|A_\mu^{-1}\,g\|_0, \]
which proves (\ref{oper6}). \hfill $\Box$ \\

 We end this Section with a number of important remarks on the theoretical and practical range of 
the splitting into regular and singular parts.
 \begin{remark}\label{rmk-wave}
Equation (\ref{laplace1}) is also useful when the {\em 2D heat equation} is considered in $\om$:
 \[\frac{\partial u}{\partial t}-\Delta u = f\mbox{ in }\om\times]0,T[,\]
 with initial condition and (homogeneous) Dirichlet boundary condition. As a matter of fact, assume
it is first discretized in time, with a time-step $\delta t$, at times $t_m=m\delta t$, $m=0,1,\cdots$:
let $u^m=u(t_m)$. Then one has to solve in space the implicit problems (with $\theta\in]0,1]$ given) \\
{\em Find $u^{m+1}\in H^1_0(\om)$ such that}
 \[
     -\Delta u^{m+1} + \frac{1}{\theta\delta t}u^{m+1} = 
 f_{(t_{m+1})}+\frac{1-\theta}{\theta}f_{(t_{m})}
       +\frac{1}{\theta\delta t}u^m+\frac{1-\theta}{\theta}\Delta u^m,
 \mbox{ in }\om.
 \]
 Above, $\theta=1$ (resp. $\theta=1/2$) corresponds to the implicit Euler (resp. Crank-Nicolson) 
 scheme. This is precisely (\ref{laplace1}) with $\mu=1/{\theta\delta t}$. \\
 Clearly, implicit schemes for the {\em 2D wave equation}
 \[\frac{\partial^2 u}{\partial t^2}-\Delta u = f\mbox{ in }\om\times]0,T[,\]
lead to other instances of Equation (\ref{laplace1}).
 \end{remark}

 \begin{remark}
Both $\phi_s$ and $p_s$ in (\ref{2ps}) are chosen independent of $\mu$, $f$ and $u_\mu$,
so their norms will be regarded as some generic constants (i.e., independent of $\mu$, $f$
and $u_\mu$.) 
\end{remark}

\begin{remark}
Instead of the decomposition (\ref{observ1}), it seems more natural 
\cite{Gris87,Gris92} to take the decomposition 
$u_\mu=\tilde u_\mu' + c_\mu\phi_\mu$, where $\phi_\mu\in H_0^1(\om)$ 
depends on the parameter $\mu$, and it is the solution to the problem: 
$-\Delta \phi_\mu + \mu\,\phi_\mu =p_\mu$ in $\om$, with 
$p_\mu\in N_\mu\setminus\{0\}$, where $N_\mu$ is given by 
\[
N_\mu=\Big\{p\in L^2(\om)\ :\ (-\Delta+\mu\,{\rm I})\,p=0,\ 
p|_{\ga_k}=0 \mbox{ in } (H_{00}^{1/2}(\ga_k))',\ 1\le k\le K\Big\}.
\]
But the decomposition (\ref{observ1}) has an important advantage: the singular part $\phi_s$ 
is independent of the parameter $\mu$. As we shall see, this will be much less expensive than
using the above more natural decomposition.
\end{remark}

 \section{Discrete formulation in the 2D domain $\om$: the SCM}
In this Section we shall formulate the generalized SCM for solving the coupled 
system~(\ref{cmu3})-(\ref{cmu1}) and derive the error estimates of the approximate solutions.
The SCM was first introduced by Assous {\it et al} \cite{AsCS98} for solving the 2D static or 
unsteady Maxwell equations without charges, and then used in \cite{CiHe03} for the 2D Poisson 
problem. As we will see, the formulation of the SCM for the 2D Poisson-like problem~(\ref{2dprob}) 
is quite different here due to the involvement of the parameter $\mu$.

Let ${\cal T}_h$ be a regular triangulation of the domain $\om$, with 
vertices $\{M_j\}_{j=1}^{N_i+N_b}$ and the last $N_b$ vertices lying 
on the boundary $\pa\om$. We define $V^h$ to be 
the continuous piecewise linear finite element space on ${\cal T}_h$ with 
the standard basis functions $\{\psi_j\}_{j=1}^{N_i+N_b}$ (cf.~\cite{Ciar91}). We further 
define $V_0^h$ to be the subspace of $V^h$ with all functions vanishing on the boundary
of $\om$. The interpolation associated with the space $V^h$ will be denoted by $\Pi_h$. 


\subsection{Approximation of the singular function $p_s$\label{sec:ps}} 
We start with the finite element approximation of the singular 
function $p_s\in N$ in (\ref{2ps}). Recall the splitting (see \cite{CiHe03})
\[
 p_s=\tilde p+p_{_P}, \quad \tilde p\in H^1(\om), 
~~p_{_P}= \rho^{-\alpha} \sin (\alpha\phi)\,.
\]
As $p_s$ is harmonic in $\om$, one can directly 
verify that the regular part $\tilde p$ in the splitting 
solves the problem:  \\
{\em Find $\tilde p\in H^1(\om)$ such that} $\tilde p = s$ on $\pa\om$ and 
\begin{equation} 
 (\nabla \tilde p, \nabla v) = 0 \quad \forall\,v\in H_0^1(\om)
\quad \label{2ps2}
\end{equation}
where the boundary function $s$ is given by 
\[
  s=0 \quad \mbox{on} \quad 
\gamma_1\cup \gamma_2; \quad s=-p_{_P}
\quad \mbox{on} \quad \gamma_k ~(3\le k\le K). 
\] 

For the finite element approximation of the problem (\ref{2ps2}), 
we shall use the simple treatment of the boundary condition: 
\begin{equation} \label{bndry1}
 \pi_h\,s=\sum_{j=N_i+1}^{N_i+N_b} s(M_j) \psi_j\,.
\end{equation}
Then we approximate $p_s$ by $p_s^h=\tilde p_h+p_{_P}$, where $\tilde p_h$ 
is the piecewise linear finite element solution to 
the problem (\ref{2ps2}). Namely,  $\tilde p_h=\pi_h s + p_h^0$ where 
$p_h^0\in V_h^0$ solves 
\begin{equation} \label{ps6}
 (\nabla \tilde p_h, \nabla v_h) = 0 \quad \forall\,v_h\in V_0^h\,.
\end{equation}

The error estimates for the singular function $p_s$ and its finite element 
approximation $p_s^h$ are summarized in the following lemma.
\begin{lemma}\label{lem:bound2} We have 
\footnote{
By construction, neither $p_s$ nor $p_s^h$ belong 
to $H^1(\om)$, due to the presence of
$p_{_P}$, but the following holds:
$$
p_s-p_s^h=\tilde p-\tilde p_h\in H^1(\om).
$$
}
\[ 
|p_s-p_s^h|_1 \lesssim h^{\alpha_0}\,,\quad
\|p_s-p_s^h\|_0 \lesssim h^{2\alpha_0}\,. 
\] 
\end{lemma}
\noindent {\it Proof}. We introduce a smooth continuation of $s$ into $\om$:
\[
 \tilde s =-p_{_P} (1-\xi(\rho))\,.
\]
Clearly, $\tilde s=s$ on $\pa\om$ and $\tilde s\in H^2(\om)$. 
Let $p^0 = \tilde p -\tilde s$. 
It is known that $\tilde p\in H^{1+\alpha_0}(\om)$, so we have 
$p^0\in H^{1+\alpha_0}(\om)\cap H_0^1(\om)$.  It follows from (\ref{2ps2}) that 
\begin{equation}\label{pint1}
 (\nabla p^0, \nabla v) =-(\nabla \tilde s, \nabla v) \quad 
 \forall\, v\in H_0^1(\om)\,.
\end{equation}
Recall $\Pi_h$ is the interpolant associated with $V^h$, thus we 
can rewrite the finite element solution 
$\tilde p_h$ to the system (\ref{ps6}) as $\tilde p_h=\Pi_h \tilde s + p_h^0$ with 
$p_h^0\in V_0^h$ now solving 
\begin{equation}\label{pint3a} 
(\nabla p^0_h,   \nabla v_h) =
  -(\nabla \Pi_h\tilde s,   \nabla v_h) \quad \forall\,v_h\in V_0^h\,,
\end{equation}
by noting $\Pi_h \tilde s= \pi_h\,s$ on $\pa\om$.

Now we are ready to derive the error estimates. 
It is clear from (\ref{pint1}) and (\ref{pint3a}) that 
\begin{equation}\label{pint3}
 (\nabla ( p^0- p^0_h),   \nabla v_h) = 
  (\nabla (\Pi_h\tilde s-\tilde s),   \nabla v_h) \quad \forall\,v_h\in V_0^h\,. 
\end{equation}
Using this, we obtain for any $q_h\in V_0^h$ that 
\[
\|\nabla (p^0-q_h)\|^2 \ge 
 \|\nabla (p^0-p_h^0)\|^2 + 2 (\nabla (\Pi_h\tilde s-\tilde s), \nabla (p^0_h-q_h)), 
\]
taking $q_h=\Pi_h p^0$ above and using the Young inequality leads to 
\begin{eqnarray*}
 |p^0-p_h^0|_1^2 &\le& |p^0-\Pi_h p^0|_1^2 + 
 2\,|\Pi_h \tilde s-\tilde s|_1\,(|p^0_h-p^0|_1+|p^0-\Pi_h p^0|_1) \\
 &\le& 2\,|p^0-\Pi_h p^0|_1^2 +\frac 12 |p^0_h-p^0|_1^2+ 
3\,|\Pi_h \tilde s-\tilde s|_1^2\,.
\end{eqnarray*}
Then by the standard interpolation results we obtain 
\[|p^0-p_h^0|_1^2\le 4\,|p^0-\Pi_h p^0|_1^2+6\,|\Pi_h \tilde s-\tilde s|_1^2
\lesssim h^{2\alpha_0}|p^0|^2_{1+\alpha_0} + h^2\,|\tilde s|^2_{2}.\]
This leads to the desired $H^1$-norm error estimate:
 \begin{eqnarray*}
 |p_s-p_s^h|_1
 &=&|\tilde p-\tilde p_h|_1=|p^0+\tilde s-p^0_h-\Pi_h\tilde s|_1\\
 &\le&|p^0-p^0_h|_1+|\tilde s-\Pi_h\tilde s|_1
  \lesssim h^{\alpha_0}+h|\tilde s|_2\lesssim h^{\alpha_0}.
 \end{eqnarray*}

Finally, we apply the Nitsche trick to derive the $L^2$-norm error estimate. 
Let $w\in H_0^1(\om)$ be the solution to the variational problem 
\begin{equation}\label{pint2}
 (\nabla w, \nabla v)=(p^0-p^0_h, v) \quad \forall\,v\in H_0^1(\om). 
\end{equation}
By the elliptic theory, we know $w\in H^{1+\alpha_0}(\om)$ and 
\[
 |w|_{1+\alpha_0} \lesssim \|p^0-p^0_h\|_0. 
\]
Let $w_h$ be the finite element approximation of 
$w$: ~$w_h\in V_0^h$ solves 
\[ 
 (\nabla w_h, \nabla v_h)=(p^0-p^0_h, v_h) \quad \forall\,v_h\in V_0^h.
\] 
Taking $v_h=w_h$ above and using the Poincar\'e inequality, 
we know 
\[
 |w_h|_1\lesssim \|p^0-p^0_h\|_0.
\]
Also, by the standard error estimate, we have 
\[
 |w-w_h|_1 \lesssim h^{\alpha_0}|w|_{1+\alpha_0}\lesssim h^{\alpha_0}\|p^0-p^0_h\|_0.
\]
Now, taking $v=p^0-p^0_h$ in (\ref{pint2}) and using (\ref{pint3}) 
and the duality argument, we obtain  
\begin{eqnarray*}
\|p^0-p^0_h\|_0^2 &=& (\nabla w, \nabla (p^0-p^0_h)) \\
 &=& (\nabla (w-w_h), \nabla (p^0-p^0_h)) + (\nabla w_h, \nabla (p^0-p^0_h))\\
 &=& (\nabla (w-w_h), \nabla (p^0-p^0_h)) 
+ (\nabla (\Pi_h \tilde s -\tilde s), \nabla (w_h-w) )\\
 &&+ (\nabla (\Pi_h \tilde s -\tilde s), \nabla w) \\
&\le& |w-w_h|_1\,|p^0-p^0_h|_1 + 
   |\Pi_h \tilde s -\tilde s|_1\,|w_h-w|_1
+ |\Pi_h \tilde s -\tilde s|_{1-\alpha_0}|w|_{1+\alpha_0}\\
&\lesssim& h^{2\alpha_0} \|p^0-p^0_h\|_0 + h^{1+\alpha_0}|\tilde s|_2\|p^0-p^0_h\|_0\,, 
\end{eqnarray*}
which leads to the desired $L^2$-norm error estimate:
 \[\hskip 20mm
 \|p_s-p_s^h\|_0
 \le\|p^0-p^0_h\|_0+\|\tilde s-\Pi_h\tilde s\|_0
  \lesssim h^{2\alpha_0}+h^2|\tilde s|_2\lesssim h^{2\alpha_0}.
 \hskip 18mm\Box\]
 \begin{remark} Following the proof given in \cite{AmMo89}, one can improve 
the results of the previous Lemma. Indeed, one can derive the estimates 
$|p_s-p_s^h|_1 \lesssim h^{\alpha}$ and $\|p_s-p_s^h\|_0 \lesssim h^{2\alpha}$, with slightly 
more restrictive assumptions on the mesh.
 \end{remark}
 \vskip -10mm
\subsection{Approximation of the singular part $\phi_s$\label{sec:phis}} 
In order to approximate the singular part $\phi_s$ in the 
decomposition $u_\mu=\tilde u_\mu+c_\mu\,\phi_s$, we recall (cf. \cite{CiHe03})
that $\phi_s\in H_0^1(\om)$ 
solves the elliptic problem (\ref{2ps}) and has the following 
decomposition: 
\begin{equation}\label{gamma1}
 \phi_s=\tilde \phi+\beta^\star\phi_{_P}, \quad \tilde \phi\in H^2(\om), 
~~\beta^\star=\frac1\pi\|p_s\|_0^2,~~\phi_{_P}= \rho^{\alpha} \sin (\alpha\phi)\,.
\end{equation}
Using (\ref{2ps}), 
we see that $\tilde \phi$, satisfying $\tilde\phi=-\beta^\star\phi_{_P}$ on 
$\partial\om$, solves the variational problem: 
\begin{equation}\label{phi1}
(\nabla \tilde \phi, \nabla v) = (p_s, v) 
 \quad \forall\, v\in H_0^1(\om).
\end{equation}

The next step is to consider the finite element approximation of $\tilde\phi$ in $V^h$:
\[\tilde\phi_h= -\beta^\star_h\pi_h\phi_P + \phi_h^0,\]
where $\pi_h$ is defined as in (\ref{bndry1}), $\beta^\star_h$ is computed using
$\beta^\star_h =\ds \frac{1}{\pi} \int_\om (p_s^h)^2 d\om$,
and $\phi_h^0\in V_0^h$ is the solution to the problem: 
\begin{equation}\label{beta3}
(\nabla\tilde\phi_h, \nabla v_h) = (p_s^h, v_h)\quad \forall\,v_h\in V_0^h.
\end{equation}

Then we propose to compute the finite element approximation of $\phi_s$ by 
\[ 
 \phi^h_s = \tilde\phi_h +\beta^\star_h\phi_{_P}\,.
\]

Below, we derive the error estimates for this approximation.
\begin{lemma}\label{lem:bound3} The following error estimates hold
\[ 
|\phi_s-\phi_s^h|_1 \lesssim h\,,\quad 
\|\phi_s-\phi_s^h\|_a \lesssim \sqrt\mu\,h\,.
\] 
\end{lemma}

\noindent {\it Proof}. 
We first estimate the error $\tilde\phi-\tilde\phi_h$. Subtracting 
(\ref{beta3}) from (\ref{phi1}) yields 
 \[ 
  (\nabla (\tilde\phi-\tilde\phi_h),  \nabla v_h) = 
 (p_s-p_s^h, v_h) \quad \forall\,v_h\in V_0^h\,, 
\] 
thus we obtain for any 
$w_h\in V^h$ satisfying $w_h-\tilde \phi_h\in V_0^h$,  
\[
 |\tilde\phi -w_h|_1^2 = |\tilde\phi -\tilde\phi_h|_1^2 + 
 |\tilde\phi_h -w_h|_1^2 + 2 (p_s-p_s^h, \tilde\phi_h -w_h),
\]
which with the Young inequality and the Poincar\'e inequality 
gives 
\begin{eqnarray}
|\tilde\phi -\tilde\phi_h|_1^2&\le& |\tilde\phi -w_h|_1^2
+2 C_P \|p_s-p_s^h\|_0( |\tilde\phi -\tilde\phi_h|_1+ |\tilde\phi -w_h|_1)
\nonumber\\
&\lesssim& 2\,|\tilde\phi -w_h|_1^2 + \frac 12 |\tilde\phi -\tilde\phi_h|_1^2 +\|p_s-p_s^h\|_0^2.  
\label{phi5}
\end{eqnarray}
Noting that $\tilde\phi=-\beta^\star\phi_{_P}$ on $\partial\om$, 
so $\beta^\star_h\Pi_h\tilde\phi=\beta^\star\tilde\phi_h$ on $\partial\om$. 
Let $w_h=\beta^\star_h\Pi_h\tilde\phi/\beta^\star$, then 
$w_h-\tilde\phi_h\in V_0^h$. With this $w_h$, we derive 
from (\ref{phi5}) and Lemma~\ref{lem:bound2} that  
\begin{eqnarray}
 |\tilde\phi -\tilde\phi_h|_1^2&\lesssim& h^{2}+ (\beta^\star)^{-2} 
 |\beta^\star\tilde\phi-\beta^\star_h\Pi_h\tilde \phi|_1^2\nonumber\\
&\lesssim& h^{2}+ |\beta^\star-\beta^\star_h|^2\,|\tilde\phi|_1^2 + 
 |\beta^\star_h|^2\,|\tilde\phi-\Pi_h\tilde\phi|_1^2.\label{eq:phi1}
\end{eqnarray}
But using the definitions of $\beta^\star$ and $\beta^\star_h$, we have
\begin{equation}\label{gamma2}
 |\beta^\star-\beta^\star_h| =\frac 1{\pi} \Big|
 \|p_s\|_0^2 -\|p_s^h\|_0^2\Big| \lesssim \|p_s-p_s^h\|_0\lesssim h^{2\alpha_0}\,. 
\end{equation}
 It follows from (\ref{eq:phi1}) and the property $\tilde\phi\in H^2(\om)$ that
\[ 
|\tilde\phi -\tilde\phi_h|_1\lesssim h\,.
\] 
This with (\ref{gamma2}) and the decompositions of $\phi_s$ and $\phi_s^h$ gives 
the desired $H^1$-norm estimate:
\begin{eqnarray*}
 |\phi_s-\phi_s^h|_1 &\le& |\tilde\phi-\tilde\phi_h|_1+ 
|\beta^\star-\beta^\star_h| \,|\phi_{_P}|_1 \lesssim h\,.
\end{eqnarray*}
 Finally, by noting that both $\phi_s$ and $\phi_s^h$ 
vanish on $\gamma_1$ and $\gamma_2$, we can apply 
the Poincar\'e inequality to the function $\phi_s-\phi_s^h$ to get 
$$
 \|\phi_s-\phi_s^h\|_0\le C_P'|\phi_s-\phi_s^h|_1\,.
$$
Then the desired estimate on $\|\phi_s-\phi_s^h\|_a$  follows from 
 \[\hskip 28mm
\|\phi_s-\phi_s^h\|_a^2=|\phi_s-\phi_s^h|_1^2+\mu\,\|\phi_s-\phi_s^h\|_0^2
  \lesssim h^2 + \mu\,h^2.
  \hskip 28mm\Box\]
 \subsection{Approximation of $\tilde u_\mu$ and $c_\mu$ in decomposition 
(\ref{observ1})}  
Noting that $\tilde u_\mu$ and $c_\mu$ solve the coupled system (\ref{cmu3}) and (\ref{cmu1}),  
it is natural 
to formulate their finite element approximations as follows: \\
{\em Find $\tilde u^h_\mu\in V_0^h$ and $c^h_\mu\in \R^1$ such that}
\begin{eqnarray}
&&a_\mu(\tilde u_\mu^h, v)
+c_\mu^h\,a_\mu(\phi_s^h, v) 
 =(f, v) \quad \forall v\in V_0^h, \label{cmu6}\\
&&\Big( \Vert p_s^h\Vert_0^2 + \mu\vert\phi_s^h\vert_1^2 \Big) c_\mu^h+
\mu\,(\tilde u_\mu^h, p_s^h)=(f, p_s^h)\,,  \label{cmu5}
\end{eqnarray}
where $\phi_s^h$ and $p_s^h$ are the finite element approximations 
of $\phi_s$ and $p_s$, see Subsect.~\ref{sec:ps}-\ref{sec:phis}.

However, this formulation requires solving a coupled system, 
and it poses some difficulty in getting the error estimates as it 
does not fall into any existing saddle-point-like framework.
Instead, we are going to propose a more efficient approximation which enables us 
to find  $\tilde u^h_\mu\in V_0^h$ and $c^h_\mu$ separately. 
In fact, we can use the formula (\ref{oper3}) to first find 
$c_\mu^h$, and then use (\ref{cmu6}) to find $\tilde u^h_\mu\in V_0^h$. 
This leads to the following algorithm to find 
$\tilde u^h_\mu\in V_0^h$ and $c^h_\mu$. Let $C^\star>0$ be a fixed constant. 

\medskip 
\underline{
{\bf SCM Algorithm} for finding $\tilde u^h_\mu\in V_0^h$ and $c^h_\mu\in \R^1$}.

\smallskip 

\underline{\bf Step 1}. {\em Find $z_\mu^h\in V_0^h$ such that}
\begin{equation}\label{zmuh}
a_\mu(z_\mu^h, v) = (f, v) \quad \forall\, v\in V_0^h\,.
\end{equation} 

\indent Compute $c_\mu^h$ as follows:
\begin{equation}
c_\mu^h =\frac{(f- \mu\,z_\mu^h, \,p_s^h)}{\Vert p_s^h\Vert_0^2}\, 
\quad \mbox{if} \quad \sqrt{\mu}< C^\star\,h^{-\frac{1}{2-\alpha_0}}\,;
\label{oper3a}
\end{equation}
\indent and 
\begin{equation}
c_\mu^h =0 \quad \mbox{if} \quad \sqrt{\mu}\ge C^\star\,h^{-\frac{1}{2-\alpha_0}}\,.
\label{oper3d}
\end{equation}

\underline{\bf Step 2}. {\em Find $\tilde u^h_\mu\in V_0^h$ such that}
\begin{equation} \label{cmu6a}
a_\mu(\tilde u_\mu^h, v)
+c_\mu^h\,a_\mu(\phi_s^h, v) 
 =(f, v) \quad \forall v\in V_0^h\,.
\end{equation}
 \begin{remark}
     In practice (see \cite{CJKL+04c}), the conditions (\ref{oper3a}-\ref{oper3d}) mean that only a 
     few coefficients $(c_\mu^h)_\mu$ need to be computed, with respect to the total number of 
     Fourier modes.
 \end{remark}
Below, we shall derive the error estimates on 
$(c_\mu-c_\mu^h)$ and $(\tilde u_\mu-\tilde u^h_\mu)$. 
Recall the formula (\ref{oper3}) for $c_\mu$: 
\begin{equation}
c_\mu =\frac{(f- \mu\,z_\mu, \,p_s)}{\Vert p_s\Vert_0^2}\,, 
\label{oper3b}
\end{equation}
where $z_\mu=A_\mu^{-1}f\in H_0^1(\om)$ solves  
\begin{equation}\label{zmu}
a_\mu(z_\mu, v) = (f, v) \quad \forall\, v\in H_0^1(\om)\,.
\end{equation} 
Clearly $z_\mu=u_\mu$, the solution to the equation (\ref{laplace1}). 
But a different notation $z_\mu$ is used here for convenience, since the numerical approximation 
$z_\mu^h$ is derived with the {\em standard} piecewise linear FEM. 

\begin{lemma}\label{lem:bound4a} For the solution $z_\mu$ to 
the problem (\ref{zmu}) and its piecewise linear finite element 
approximation $z_\mu^h$ in (\ref{zmuh}), 
we have the following error estimates 
\begin{eqnarray}
\|z_\mu-z_\mu^h\|_0&\le& \mu^{-1}\|f\|_0\,, \label{zmu5}\\ 
\|z_\mu-z_\mu^h\|_0&\lesssim& h^{2\alpha_0}\mu^{\alpha_0-1}(1+\sqrt{\mu}h)^2 \|f\|_0\,, \label{zmu4}
\end{eqnarray}
while for the coefficient $c_\mu$ in (\ref{oper3b}) and its 
approximation $c_\mu^h$ in (\ref{oper3a}), we have 
\begin{equation}\label{error1} 
|c_\mu-c_\mu^h|
 \lesssim (h^{2\alpha_0}\mu^{\alpha_0}(1+\sqrt{\mu}h)^2  + h)\,\|f\|_0 \,.
\end{equation}
\end{lemma}

\noindent {\it Proof}. It follows from (\ref{zmu}) and (\ref{zmuh}) 
that 
\[
 a_\mu(z_\mu-z_\mu^h, z_\mu-z_\mu^h)=a_\mu(z_\mu, z_\mu-z_\mu^h) 
 =(f, z_\mu-z_\mu^h).
\]
This implies 
\[
 |z_\mu-z_\mu^h|^2_1+ \mu\, \|z_\mu-z_\mu^h\|^2_0 
  \le \|f\|_0\,\|z_\mu-z_\mu^h\|_0, 
\]
thus (\ref{zmu5}) follows by the Young inequality. 

We next show (\ref{zmu4}). Again 
it follows from (\ref{zmuh}) and (\ref{zmu}) that 
\[
 \|z_\mu-z_\mu^h\|_a\le \|z_\mu - v_h\|_a \quad 
 \forall\, v_h\in V_0^h. 
\]
 But, by standard interpolation theory, we know that
 \[|z_\mu-\Pi_h z_\mu|_1 \lesssim h^{\alpha_0} |z_\mu|_{1+\alpha_0},\mbox{ and }
  \|z_\mu-\Pi_h z_\mu\|_0\lesssim h^{1+\alpha_0} |z_\mu|_{1+\alpha_0}.\]
 Therefore, we reach
\begin{equation}
\|z_\mu-z_\mu^h\|_a \lesssim (1+ \sqrt{\mu}\,h) h^{\alpha_0}
 |z_\mu|_{1+\alpha_0}\,. 
\label{zmu2}
\end{equation}

 On the other hand, for any $g\in L^2(\om)$, define $w\in H_0^1(\om)$ such that 
\begin{equation}
 a_\mu(w, v)=(g, v) \quad \forall\, v\in H_0^1(\om). 
\label{zmu3}
\end{equation}
Using the duality and (\ref{zmu3}), we have 
\begin{eqnarray*}
 \|z_\mu-z_\mu^h\|_0 &=&\sup_{g\in L^2(\om)} \frac{(z_\mu-z_\mu^h, g)}{\|g\|_0}
                      = \sup_{g\in L^2(\om)} \frac{a_\mu (w, z_\mu-z_\mu^h)}{\|g\|_0}\\
                &=& \sup_{g\in L^2(\om)} \frac{a_\mu (w-\Pi_h w, z_\mu-z_\mu^h)}{\|g\|_0}
 \le \sup_{g\in L^2(\om)} \frac{\|w-\Pi_h w\|_a\,\|z_\mu-z_\mu^h\|_a}{\|g\|_0}.
\end{eqnarray*}
Using the interpolation result and the same derivation as 
in (\ref{zmu2}) and the {\em a priori} estimate (\ref{asympu}) 
(with $u$ and $f$ replaced by $w$ and $g$), we obtain  
\begin{eqnarray*}
\|z_\mu-z_\mu^h\|_0 
&\lesssim& \sup_{g\in L^2(\om)} \frac{h^{2\alpha_0}(1+\sqrt{\mu}h)^2 |w|_{1+\alpha_0}
|z_\mu|_{1+\alpha_0}}{\|g\|_0}\\
 &\lesssim&h^{2\alpha_0}\mu{^{\alpha_0-1}}(1+\sqrt{\mu}h)^2 \|f\|_0\,,
\end{eqnarray*}
which proves (\ref{zmu4}).

It remains to prove (\ref{error1}). 
We have from (\ref{oper3a}) and (\ref{oper3b}) that 
\begin{eqnarray*}
c_\mu -c_\mu^h &=& 
 \frac{(f- \mu\,z_\mu, \,p_s)}{\Vert p_s\Vert_0^2}-
 \frac{(f- \mu\,z_\mu^h, \,p_s^h)}{\Vert p_s^h\Vert_0^2}\\
&=& \Big\{ \frac{(f, \,p_s)}{\Vert p_s\Vert_0^2} - 
   \frac{(f, \,p_s^h)}{\Vert p_s^h\Vert_0^2}\Big\} +\mu\, 
  \Big\{ \frac{(z_\mu^h, \,p_s^h)}{\Vert p_s^h\Vert_0^2}-
 \frac{(z_\mu, \,p_s)}{\Vert p_s\Vert_0^2}\Big\}
: = I_1 + I_2\,.
\end{eqnarray*}
For $I_1$, we have from Lemma\,\ref{lem:bound2} that 
\[
 |I_1|\lesssim h\,\|f\|_0\,.
\]
For $I_2$, we further write it as follows 
\[
 I_2 = \mu\, \frac{(z_\mu^h-z_\mu, \,p_s^h)}{\Vert p_s^h\Vert_0^2}
   + \mu\, \frac{(z_\mu, \,p_s^h-p_s)}{\Vert p_s^h\Vert_0^2} + 
\mu\,(z_\mu, \,p_s)\Big\{\frac{1}{\Vert p_s^h\Vert_0^2} 
 - \frac{1}{\Vert p_s\Vert_0^2}\Big\}\,.
\]
Then using estimate (\ref{zmu4}) and Lemma~\ref{lem:bound2}, we can derive 
\[
 |I_2| \lesssim  \mu\,\|z_\mu-z_\mu^h\|_0 + h\,\|f\|_0 
 \lesssim (h^{2\alpha_0}\mu^{\alpha_0}(1+\sqrt{\mu}h)^2 + h)\,\|f\|_0 \,.
\]
This with the estimate of $I_1$ gives (\ref{error1}).
 \hfill $\Box$
  
\medskip
In the rest of this Section, we shall estimate the error between 
the solution $u_\mu$ to the elliptic problem (\ref{laplace1}) and its 
SCM approximation $u_\mu^h$. We note that the decomposition of $u_\mu$ is equal to:
\begin{equation}\label{u1}
u_\mu=\tilde u_\mu +c_\mu\,\phi_s= \tilde u_\mu +c_\mu\, 
 (\tilde \phi+\beta^\star \phi_{_P})\,.
\end{equation} 
So, we propose its SCM approximation $u_\mu^h$ of the form: 
\begin{equation}\label{uh1}
u_\mu^h=\tilde u_\mu^h +c_\mu^h\,\phi_s^h 
 = \tilde u_\mu^h +c_\mu^h\,(\tilde \phi_h+\beta^\star_h \phi_{_P}).  
\end{equation}

We shall derive the error estimate on $u_\mu-u_\mu^h$. Let us start 
with the estimate of $(\tilde u_\mu-\tilde u_\mu^h)$. We have 
\begin{lemma}\label{lem:bound4} The following error estimate holds
\[ 
\|\tilde u_\mu-\tilde u_\mu^h\|_a^2 \lesssim \sqrt{\mu}\,
(h^2\,\|f\|_0^2+|c_\mu-c_\mu^h|^2)\,. 
\] 
\end{lemma}

\noindent {\it Proof}. 
Subtracting (\ref{cmu6}) from (\ref{cmu3})  we have 
\[
 a_\mu(\tilde u_\mu-\tilde u_\mu^h, v_h) +c_\mu a_\mu(\phi_s, v_h) 
  -c_\mu^h a_\mu(\phi_s^h, v_h)=0 \quad \forall\, v_h\in V_0^h.
\]
Using this we obtain for any $w_h\in V_0^h$,
\[
\|\tilde u_\mu-w_h\|_a^2 = \|\tilde u_\mu-\tilde u_\mu^h\|_a^2 + 
\|\tilde u_\mu^h-w_h\|_a^2 + 2 c_\mu^h a_\mu(\phi_s^h, \tilde u_\mu^h-w_h) -2 
c_\mu a_\mu(\phi_s, \tilde u_\mu^h-w_h), 
\] 
which implies 
\begin{eqnarray}
&&\|\tilde u_\mu-\tilde u_\mu^h\|_a^2 \nonumber\\
&\le& \|\tilde u_\mu-w_h\|_a^2+ 2 c_\mu
 \,a_\mu(\phi_s-\phi_s^h, \tilde u_\mu^h-w_h) + 
 2(c_\mu-c_\mu^h) a_\mu(\phi_s^h, \tilde u_\mu^h-w_h)\label{zmu8}\\
&\le&  \|\tilde u_\mu-w_h\|_a^2+ 2 \,|c_\mu|\,
 \|\phi_s-\phi_s^h\|_a\,\|\tilde u_\mu^h-w_h\|_a + 
 2 \,|c_\mu-c_\mu^h|\, \|\phi_s^h\|_a\,\|\tilde u_\mu^h-w_h\|_a\,.\nonumber
\end{eqnarray}
Now, there holds
$\|\phi_s\|_a - \|\phi_s-\phi_s^h\|_a\le \|\phi_s^h\|_a\le \|\phi_s\|_a + \|\phi_s-\phi_s^h\|_a$.
Using Lemma~\ref{lem:bound3} and 
$\|\phi_s\|^2_a=|\phi_s|_1^2+\mu\,\|\phi_s\|_0^2$, we find
$\|\phi_s^h\|_a\approx\sqrt\mu\,\|\phi_s\|_0$.
Using the interpolation results, we obtain
\[
 \|\tilde u_\mu-\Pi_h \tilde u_\mu\|_a^2 \le 
 |\tilde u_\mu-\Pi_h \tilde u_\mu|_1^2 +
 \mu\,\|\tilde u_\mu-\Pi_h \tilde u_\mu\|_0^2
 \lesssim h^2\,|\tilde u_\mu|_2^2, 
\]
thus letting $w_h=\Pi_h \tilde u_\mu$ in (\ref{zmu8}) and using 
Lemma~\ref{lem:bound1}, we derive 
\begin{eqnarray*}
\hskip 18mm\|\tilde u_\mu-\tilde u_\mu^h\|_a^2
&\lesssim& h^2\,|\tilde u_\mu|_2^2 +  
\sqrt{\mu} h^2 \|f\|_0 \,|\tilde u_\mu|_2 + 
\sqrt{\mu}\,h\,|c_\mu-c_\mu^h|\,|\tilde u_\mu|_2 \\ 
&\lesssim& \sqrt{\mu}\,(h^2\,\|f\|_0^2+ |c_\mu-c_\mu^h|^2)\,.
\hskip 52mm\Box
\end{eqnarray*}

\begin{theorem}\label{thm:main1}
Let $u_\mu$ be the solution to the equation (\ref{laplace1})
and $u_\mu^h$ be its finite element approximation given in (\ref{uh1}). 
Then the following error estimate holds:
\[ 
\exists C>0\mbox{ such that }\forall\mu,\ \|u_\mu-u_\mu^h\|_a \le C\,\mu\,h\,\|f\|_0.
\] 
\end{theorem}

\noindent {\it Proof}. 
It follows from (\ref{u1}) and (\ref{uh1}) that 
\begin{eqnarray*}
u_\mu-u_\mu^h = (\tilde u_\mu-\tilde u_\mu^h) +c_\mu (\phi_s- \phi^h_s)
  + \phi_s^h (c_\mu-c_\mu^h). 
\end{eqnarray*}
Then we obtain, using Lemmas~\ref{lem:bound4}, \ref{lem:bound3} and \ref{lem:bound1}, that 
\begin{eqnarray*}
\|u_\mu-u_\mu^h\|^2_a &\le& 3\Big\{ \|\tilde u_\mu-\tilde u_\mu^h\|^2_a +
 |c_\mu|^2\,\| \phi_s- \phi_s^h\|_a^2
  + \|\phi_s^h\|_a^2 |c_\mu-c_\mu^h|^2\Big\} \\ 
 \label{uuh2}
 &\lesssim&  \mu\,h^2\,\|f\|_0^2 +\mu\,|c_\mu-c_\mu^h|^2\,. 
\end{eqnarray*}
To prove the desired estimate, we need simply
\begin{equation}
 |c_\mu-c_\mu^h|^2 \lesssim \mu\,h^2\,\|f\|_0^2\,. \label{eq:err4}
\end{equation}

First consider the case (\ref{oper3d}), i.e., $\sqrt{\mu}\ge C^\star\,h^{-\frac{1}{2-\alpha_0}}$.
This condition is equivalent to
\[h^{-2}\mu^{\alpha_0-2}\lesssim 1.\]
Then (\ref{eq:err4}) comes directly from this condition, $c_\mu^h=0$ and (\ref{asympc}) as follows:
\[
 |c_\mu-c_\mu^h|^2 = c_\mu^2 \lesssim \mu^{\alpha_0-1}\|f\|_0^2 
\lesssim {\mu}\,h^2\,(h^{-2}\mu^{\alpha_0-2}) \|f\|_0^2 
 \lesssim {\mu}\,h^2\,\|f\|_0^2\,. 
\]

For the remaining case (\ref{oper3a}), we have $\sqrt{\mu}<  C^\star\,h^{-\frac{1}{2-\alpha_0}}$,
or $h^2\lesssim \mu^{-(2-\alpha_0)}$. 
On the one hand, since $\alpha_0<1$, $\sqrt{\mu}h\lesssim h^{\frac{1-\alpha_0}{2-\alpha_0}}\lesssim 1$.
On the other hand, since $2\alpha_0-1>0$, 
$h^{4\alpha_0-2}\lesssim \mu^{-(2\alpha_0-1)(2-\alpha_0)}$. But one infers from (\ref{error1})
and these inequalities that 
\begin{eqnarray*}
|c_\mu-c_\mu^h|^2 
 &\lesssim &(h^{4\alpha_0}\mu^{2\alpha_0} + h^2)\,\|f\|^2_0
 \lesssim h^2\,(\mu^{2\alpha_0-(2\alpha_0-1)(2-\alpha_0)}+1)\,\|f\|^2_0\,.
\end{eqnarray*}
To conclude, (\ref{eq:err4}) follows from this and the fact that, as $\alpha_0\in]\frac12,1[$,
the exponent of $\mu$ is bounded by 
\[\hskip 12mm
 2\alpha_0-(2\alpha_0-1)(2-\alpha_0) = 2\alpha_0^2-3\alpha_0+ 2 
 =1+ (2\alpha_0-1)(\alpha_0-1) < 1\,.
 \hskip 12mm\Box\]

\section{Fourier Singular Complement Methods}\label{FSCM}
 In order to define the numerical part of the Fourier Singular Complement Method,
let us prove a result which can be viewed as the mathematical foundation of the FSCM, 
from the Fourier point of view. It allows to recover (\ref{prismatic-splitting}-\ref{reg-u-2}), 
for sufficiently smooth right-hand sides. \\

 Let $u$ be the solution to the Poisson problem (\ref{Poisson}) and $u_k$ be its 
Fourier coefficients in (\ref{fourier0}). By Lemma~\ref{lem:f4}, we know that 
$u_k(x_1, x_2)$ solves the 2D problem (\ref{Poissonk}-\ref{2dprob}). And using 
(\ref{observ1}) we can decompose $u_k$ as follows:
\begin{equation}\label{uk-decomp}
 u_k=\tilde u_k + c_k\, \phi_s 
\end{equation}
where $\tilde u_k\in H^2(\om)\cap H_0^1(\om)$ and $\phi_s\in H_0^1(\om)$ 
solves (\ref{2ps}). 

\begin{lemma}\label{lem:f6}
Let $f\in h^2(\Om)\cap h_{\diamond}^1(\Om)$, and $u\in H_0^1(\Om)$ be the solution to 
(\ref{Poisson}). Then 
 \begin{equation}\label{prismatic-splitting2}
     u=\tilde u + \ga(x_3)\phi_s,\mbox{ with }
     \tilde u\in H^2(\Om)\cap H^1_0(\Om),\ 
     \ga\in H^2(]0,L[)\cap H^1_0(]0,L[).
 \end{equation}
\end{lemma}
  \noindent {\it Proof}. Let $(U_K)_K$ be the Fourier sequence of $u$. Recall that $(U_K)_K$ 
converges to $u$ in $H^1_0(\Om)$, and $(\Delta U_K)_K$ converges to $-f$ in $L^2(\Om)$. From 
(\ref{uk-decomp}), let us split the Fourier sequence into regular and singular parts, as
 \[U_K = \widetilde U_K + \ga_K(x_3)\,\phi_s ,\mbox{ with }
         \widetilde U_K = \sum_{k=1}^K \tilde u_k\sin\frac{k\pi}{L}x_3,\ 
             \ga_K(x_3) = \sum_{k=1}^K c_k\sin\frac{k\pi}{L}x_3. \]
 We shall prove below that $(\ga_K)_K$ converges in $H^2(]0,L[)\cap H^1_0(]0,L[)$, and 
$(\widetilde U_K)_K$ converges in $H^2(\Om)\cap H^1_0(\Om)$.  \\

 As far as the singular part is concerned, from (\ref{h10-cap-h2}) and the bound on $|c_k|$ in 
Lemma~\ref{lem:bound1}, we obtain that $\ds\sum_{k=1}^\infty k^4 |c_k|^2<\infty$. Since we are 
dealing with the 1D Fourier sequence $(\ga_K)_K$ (with sine functions), it is well-known that 
it converges to a limit, subsequently called $\ga$, in $H^2(]0,L[)\cap H^1_0(]0,L[)$. Then, 
one finds that $(\ga_K\,\phi_s)_K$ converges to $\ga\,\phi_s$ in $H^1_0(\Om)$, and that 
$(\Delta(\ga_K\,\phi_s))_K$ converges in $L^2(\Om)$, to $\ga''\,\phi_s-\ga p_s$. \\

For the regular part, we note that since there holds $\widetilde U_K=U_K-\ga_K\,\phi_s$, 
$(\widetilde U_K)_K$ converges in $H^1_0(\Om)$, to a limit called $\tilde u$, which is equal to
\[\tilde u=u-\ga\,\phi_s.\] 
Moreover, $(\Delta\widetilde U_K)_K$ converges in $L^2(\Om)$, to $\Delta\tilde u$. \\
 
 To conclude the proof, one has to establish that $\tilde u$ is an element of $H^2(\Om)$. From 
Corollary \ref{cor:f2}, we know already that $\pa_3\tilde u$ is in $H^1(\Om)$. So one has to check 
that $\pa_{ij}\tilde u$ is in $L^2(\Om)$, for $i,j\in\{1,2\}$. But this follows from the estimate 
on $|\tilde u_k|_2$ in Lemma~\ref{lem:bound1}, and on the expression of the second order partial 
derivatives of $\widetilde U_K$, that is
 \[\hskip 48mm
 \pa_{ij}\widetilde U_K = \sum_{k=1}^K \pa_{ij}\tilde u_k\sin\frac{k\pi}{L}x_3.
 \hskip 48mm\Box \]
 \begin{remark} In the more general case, i.e., $f\in L^2(\Om)$, one gets only a convergence of 
$(\ga_K)_K$ in $H^{1-\alpha}(]0,L[)$, see \cite{BrNS0x}. This precludes a convergence of the 
singular part in the desired Sobolev spaces, i.e., $H^1(\Om)$ with $L^2(\Om)$ Laplacian.
 \end{remark}
In order to build the numerical schemes which completely define the FSCM, we introduce 
$u_{k}^h(x_1, x_2)$ the SCM approximation to $u_k(x_1, x_2)$. It is the same as 
$u_\mu^h$ in (\ref{uh1}), but with $\mu$ replaced by $k^2\pi^2/L^2$, that is, 
\[
 u_k^h =\tilde u_k^h + c_k^h\, \phi_s^h.
\]

We then rephrase the 2D SCM Algorithm (\ref{zmuh}-\ref{cmu6a}). This gives \\

\underline{\bf Step 1}. {\em Find $z_k^h\in V_0^h$ such that}
\begin{equation}\label{zmuh-bis}
a_k(z_k^h, v) = (f, v) \quad \forall\, v\in V_0^h\,.
\end{equation} 

\indent Compute $c_k^h$ as follows:
\begin{equation}
c_k^h =\frac{(f- \ds\frac{k^2\pi^2}{L^2}\,z_k^h, \,p_s^h)}{\Vert p_s^h\Vert_0^2}\, 
\quad \mbox{if} \quad k< C^\star\frac{L}{\pi}\,h^{-\frac{1}{2-\alpha_0}}\,;
\label{oper3a-bis}
\end{equation}
\indent and 
\begin{equation}
c_k^h =0 \quad \mbox{if} \quad k\ge C^\star\frac{L}{\pi}\,h^{-\frac{1}{2-\alpha_0}}\,.
\label{oper3d-bis}
\end{equation}

\underline{\bf Step 2}. {\em Find $\tilde u^h_k\in V_0^h$ such that}
\begin{equation} \label{cmu6a-bis}
a_k(\tilde u_k^h, v)
+c_k^h\,a_k(\phi_s^h, v) 
 =(f, v) \quad \forall v\in V_0^h\,.
\end{equation}

 As mentioned already, only a few coefficients $(c_k^h)_k$ are actually computed. \\
 
 Following (\ref{FSCM-approx}), we finally define the FSCM approximation to the solution $u$ 
to (\ref{Poisson}) as follows: 
\[ 
 U_N^h (x_1, x_2, x_3) =\sum_{k=1}^N u_{k}^h(x_1, x_2) 
 \sin \frac{k\pi}{L}x_3\,.
\] 
Then we have the final error estimate below
\begin{theorem} 
Assume that $f\in h^1_\diamond(\Om)\cap h^2(\Om)$. \\
The following error estimate holds: 
\[
 \|\nabla (u-U_N^h)\|_{L^2(\Om)} 
\lesssim (h+N^{-1})\Big\{ \|f\|_{L^2(\Om)} +
\|\pa_{33}f\|_{L^2(\Om)} \Big\}\,.
\]
\end{theorem}
\noindent {\it Proof}. 
Using the Fourier expansion of $u$ and the definition of $U_N^h$, we have, cf. (\ref{f-H1_norm}),
\begin{eqnarray*}
 \|\nabla (u-U_N^h)\|_{L^2(\Om)}^2 
 &=& \frac{L}2 \sum_{k=1}^N \Big(\|\nabla ( u_k-u_{k}^h)\|_0^2 +
 (\frac{k\pi}{L})^2 \|u_k-u_{k}^h\|_0^2 \Big)\\
 &&+ \frac{L}2 \sum_{k>N} \Big(\|\nabla  u_k\|_0^2 +
 (\frac{k\pi}{L})^2 \|u_k\|_0^2 \Big)\\
 &=:& \mbox{I}_1 +\mbox{I}_2.
\end{eqnarray*}
According to Lemma~\ref{lem:f4}, we derive  
\begin{eqnarray*} 
\mbox{I}_2 &=& \frac{L}2 \sum_{k>N} \Big(\|\nabla  u_k\|_0^2 +
 (\frac{k\pi}{L})^2 \|u_k\|_0^2 \Big)\\
 &\le& \frac{L}2 N^{-2} \sum_{k>N} k^2 \Big(\|\nabla  u_k\|_0^2 +
 (\frac{k\pi}{L})^2 \|u_k\|_0^2 \Big)\\
&\le& \Big( \frac{L}{\pi}\Big)^2 N^{-2}\,\|f\|_{L^2(\Om)}^2\,. 
\end{eqnarray*} 
For $\mbox{I}_1$, we have 
\[
 \mbox{I}_1= \frac{L}2 \sum_{k=1}^N \|u_k-u_{k}^h\|^2_a\,. 
\]
According to Theorem~\ref{thm:main1} we have 
\[
 \|u_k-u_{k}^h\|^2_a \lesssim k^4\, h^2\,\|f_k\|_0^2\,.
\]
Using this and (\ref{h10-cap-h2}), we obtain the estimate of $\mbox{I}_1$: 
\[
 \mbox{I}_1\lesssim h^2\,\sum_{k=1}^N k^4\,\|f_k\|^2_0 
 \lesssim h^2\,\|\pa_{33}f\|_{L^2(\Om)}^2,
\]
which, together with the previous estimate of $\mbox{I}_2$, leads to 
the desired error estimate. 
 \hfill $\Box$

\section{Conclusion}\label{Conclusion}
 The optimal convergence rate of the FSCM in prismatic domains, has been proven for the Poisson 
problem with homogeneous Dirichlet boundary conditions. Assuming that the right-hand side $f$ is
slightly more regular than  $f\in L^2(\Om)$, i.e., that $f$ belongs to 
$h^2(\Om)\cap h_{\diamond}^1(\Om)$, the convergence rate of the FSCM in $H^1$-norm
is like
 \[\|u-U_N^h\|_1\le C_f (h+N^{-1}),\]
 where $h$ is the 2D mesh size, and $N$ is the number of Fourier modes used. \\
 
 The same result also holds for the discretization of the Poisson problem with a homogeneous
Neumann boundary condition, or for the Poisson problem with non-homogeneous boundary conditions, 
provided there exist sufficiently smooth liftings. \\

 Further, it is no difficulty to consider the case of a prismatic domain $\Om$ with several
reentrant edges, i.e., $\om$ with several reentrant corners. \\
 
 As far as the assumptions on the right-hand side $f$ are concerned, a few remarks can be made. 
 It seems that, in a prismatic domain $\Om$ such as the one we considered here, the boundary condition 
on the bases was omitted in \cite{Apel99}. 
 Nevertheless, this condition does not exist in the case of an axisymmetric domain, see 
\cite{CJKL+04b}, nor in the case of an infinite cylinder. In other words, $f\in h^2(\Om)$ is enough 
in those types of domains. 
 In the case of a Poisson problem with Neumann boundary conditions, one has to replace the vanishing 
trace conditions at the bases by the familiar $\pa_3 f=0$ at the same bases. \\

 As mentioned already, this paper is the first part of a three-part article \cite{CJKL+04b,CJKL+04c}.
In the companion paper \cite{CJKL+04b}, the FSCM is analysed theoretically and its numerical 
approximation is built, in {\em axisymmetric domains with conical vertices and reentrant edges}. There 
are two difficulties which are inherent in this class of domains. The first one is the weights, 
which have to be introduced in the 2D sections. The second one is the addition of {\em sharp vertex} 
singularities, which have to be taken into account separately. 
 In \cite{CJKL+04c}, the FSCM is analyzed from a numerical point of view (complexity, implementation 
issues, numerical experiments, etc.), and it is compared to other methods, such as mesh refinement 
techniques, or variants of the FSCM (2D SCM with the $\lambda$-approach \cite{CiHe03}; 3D discretization 
of the regular part, etc.) in  prismatic or axisymmetric domains. In particular, the use of the FFT 
to aproximate the sine functions in $x_3$ is motivated and justified there. \\

 As noted in Remark \ref{rmk-wave}, one can apply the same theoretical and numerical techniques to the 2D 
heat or wave equations, with any $L^2$-smooth (in space) right-hand side. For these PDEs, the 
singular functions $p_s$ and $\phi_s$ do not depend on the time-step. \\
Finally, the results, can also be viewed as the first effort towards the discretization of electromagnetic 
fields in prismatic domains, with {\em continuous numerical approximations}, the importance of which is well-known,
cf. \cite{ADHR93}. As a matter of fact, the SCM developed in \cite{AsCS98,AsCS00,Garc02} for 2D 
electromagnetic computations can be generalized, based on the results obtained here.

\end{document}